\documentclass{article}
\usepackage{amsmath,amsfonts,epsfig}
\usepackage{mathptmx}
\usepackage{amsmath,amsfonts,latexsym,amssymb}
\usepackage{mathrsfs}
\usepackage{verbatim}
\usepackage[utf8]{inputenc} 
\usepackage[T1]{fontenc}
\usepackage{ae,aecompl}
\usepackage{braket}

\newcommand{\bb}{\mathbb}

\newcommand{\cx}{{\bb C}}

\renewcommand{\bold}[1]{\medskip \noindent {\bf \boldmath #1
                        }\nopagebreak[4]}

\newcommand{\del}{\partial}

\newcommand{\chat}{\widehat{\cx}}

\newcommand{\area}{\operatorname{area}}

\newcommand{\inj}{\operatorname{inj}}
\newcommand{\id}{\operatorname{id}}

\renewcommand{\Im}{\operatorname{Im}}

\renewcommand{\Re}{\operatorname{Re}}

\newcommand{\sech}{\operatorname{sech}}

\newcommand{\supp}{\operatorname{supp}}

\newcommand{\Teich}{\operatorname{Teich}}

\newcommand{\vol}{\operatorname{vol}}

\newtheorem{theorem}{Theorem}[section]
\newtheorem{prop}[theorem]{Proposition}
\newtheorem{lemma}[theorem]{Lemma}

\newcommand{\cA}{{\cal A}}

\newcommand{\cC}{{\cal C}}
\newcommand{\cD}{{\cal D}}

\newcommand{\cF}{{\cal F}}

\newcommand{\cL}{{\cal L}}

\newcommand{\cO}{{\cal O}}

\newcommand{\cS}{{\cal S}}

\newcommand{\cU}{{\cal U}}

\newcommand{\calC}{{\mathcal C}}

\newcommand{\calS}{{\mathcal S}}

\newtheorem{corollary}[theorem]{Corollary}

\renewcommand{\hbar}{\bar{{\mathbb H}}^3}

\newcommand{\CC}{\mathbb C}

\newcommand{\R}{\mathbb R}

\newcommand{\Z}{\mathbb Z}

\newcommand{\Hp}{{\mathbb H}^2}
\newcommand{\Hs}{{{\mathbb H}^3}}

\newcommand{\Rvol}{{{\rm Vol}_R}}
\newcommand{\Cvol}{{{\rm Vol}_C}}

\newcommand{\psl}{{\rm PSL}(2,\CC)}

\def\eproof{$\Box$ \medskip}
\renewcommand{\area}{{\operatorname{{\bf area}}}}

\newcommand{\bx}{\mathcal B_X}

\newcommand{\bC}{{\bf C}}
\newcommand{\bS}{{\bf S}}
\newcommand{\bT}{{\bf T}}
\newcommand{\barHs}{{\bar{\mathbb H}^3}}

\begin{document}

\title{\bf Convergence of the   gradient flow of renormalized volume to convex cores with totally geodesic boundary}
\author{Martin
   Bridgeman, \   Kenneth Bromberg, \ and
  Franco Vargas Pallete\thanks{M. Bridgeman's research  was supported by NSF grants  DMS-1564410 and DMS-2005498. K. Bromberg's research was supported by NSF grants
    DMS-150917 and, DMS-1906095. F. Vargas Pallete's research was supported by NSF grant DMS-2001997. This paper was also supported by the National Science Foundation under Grant No. 1440140, DMS-1928930, while  Bridgeman and Vargas Pallete were visiting researchers at the Mathematical Sciences Research Institute in Berkeley, California, during Fall 2020.}}

\date{\today}

\maketitle
\begin{abstract}
We consider   the Weil-Petersson gradient vector field of renormalized volume  on the deformation space of convex cocompact hyperbolic structures on (relatively) acylindrical manifolds. In this paper we prove the conjecture that the flow has a global attracting fixed point  at the structure $M_{\rm geod}$ the unique structure with minimum convex core volume.
\end{abstract}

\section{Introduction}
The deformation space of convex cocompact structures $CC(N)$ on a hyperbolizable 3-manifold $N$ has a natural flow $V$, first studied in \cite{BBB}. This flow $V$ has a classical description; at a point $M \in CC(N)$ it is the Weil-Petersson dual  of the Schwarzian derivative of the maps uniformizing the components of the conformal boundary $\partial_cM$ of $M$. Work of Storm proved that the convex core volume is minimized if and only if $N$ is acylindrical with the minimum given by the manifold $M_{\rm geod}$ whose convex core boundary is totally geodesic (see \cite{storm1}). A natural conjecture is that the flow  $V$ uniformizes $N$. Specifically that for any flowline $M_t$ of $V$ we have $M_t \rightarrow M_{\rm geod}$. In this paper, we prove this conjecture and extend it to the class of relatively acylindrical manifolds.

Although the flow has the above classical description in terms of the Schwarzian derivative, it only arose recently in the study of renormalized volume. This perspective will not be needed in this paper, but renormalized volume gives an analytic function  ${\Rvol}:CC(N)\rightarrow \R$ and the flow $V$ is equal to the Weil-Petersson gradient flow of $-{\Rvol}$. 
Renormalized volume was  introduced in work of Graham and Witten (\cite{GW}) in physics to give an
alternative notion of volume for conformally compact Einstein manifolds. In the hyperbolic setting, this was described and developed in the papers 
\cite{TT,ZT,KS08,KS12}) of Takhtajan, Zograf, Teo, Krasnov, and
Schlenker. The renormalized volume $\Rvol(M)$ of a hyperbolic manifold $M$ connects
many analytic concepts from the deformation theory with the geometry of
$M$ and is closely related to classical objects such as the convex
core volume $\Cvol(M)$ and the Weil-Petersson geometry of  Teichm\"uller
space. For the description of these connections, we refer the reader to the earlier papers \cite{BBB,wp-paper} for this perspective.

\subsection{Flow on a deformation space of relatively acylindrical manifold}

For $N$ a compact hyperbolizable 3-manifold, we denote by $CC(N)$ the space of convex cocompact hyperbolic structures on the interior of $N$. We consider triples $(N;S,X)$ where $S$ is a union of components of $\partial N$ and $X$  a conformal structure on $\partial N-S$. We then define $CC(N;S,X) \subseteq CC(N)$ to be the subset with conformal structure $X$ on $\partial N-S$.  The pair $(N;S)$ is {\em relatively acylindrical} if there are no non-trivial annuli with boundary curves both in $S$. There are two important examples. The first is when $S = \del N$ and $N$ is a acylindrical. The second important example is the pair $(S\times[0,1],S\times\{0\})$ where $S$ is a closed surface. Then $CC(S\times[0,1];S\times\{0\},X)$ is called the {\em Bers slice} and denoted $\bx$. While these may be the two main cases of interest, our result will hold in the general setting of $CC(N;S,X)$ for any relatively acylindrical $(N;S)$. By the classical deformation theory of Kleinian groups (see \cite{Kra:crash}), $CC(N;S,X)$ is parameterized by the Teichm\"uller space $\Teich(S)$ of conformal or hyperbolic structures on $S$. Combined with Thurston's hyperbolization theorem this deformation theory also implies that 
 there is a unique $M_{\rm geod}   \in CC(N;S,X)$ such that $\partial_c M_{\rm geod} = X \sqcup Y_{\rm geod}$ and the convex core  of $M_{\rm geod}$ has totally geodesic boundary facing $Y_{\rm geod}$.

Given $Y\in\Teich(S)$ we let  $M _Y\in CC(N;S,X)$ be the convex cocompact hyperbolic 3-manifold whose conformal boundary restricted to $S$ is $Y$. Let $\Gamma_Y$ be a Kleinian group with $M_Y = \Hs/\Gamma_Y$ and let $\Omega_Y$ be the union of components of the domain of discontinuity of $\Gamma_Y$ that project to $Y$. The components of $\Omega_Y$ will be Jordan domains in $\chat$. Given a component $\Omega$ of $\Omega_Y$ let $f\colon \Hp\to \Omega$ be a uniformizing univalent map. The {\em Schwarzian derivative} $S(f)$ defines a holomorphic quadratic differential on $\Omega$. If we repeat this construction for every component of $\Omega_Y$ we get a $\Gamma_Y$-invariant holomorphic quadratic on $\Omega$ which will descend to a holomorphic quadratic differential $\phi_Y$ on $Y$.

Recalling that we have an isomorphism $CC(N;S,X) \cong \Teich(S)$ and that tangent vectors in $T_Y\Teich(S)$ are given by Beltrami differentials on $Y$ we can define a vector field $V= V_{(N;S,X)}$ on $\Teich(S)$ by taking the harmonic Beltrami differential associated to $\phi_Y$. Namely let 
$$V(Y) = -\left[ \frac{\bar\phi_Y}{\rho_Y}\right]$$
where $\rho_Y$ is the area form for the hyperbolic metric on $Y$. The expression inside the brackets is a Beltrami differential with the brackets indicating that we are taking the equivalence class in the tangent space $T_Y\Teich(S)$. Thus $V$ is a vector field on $\Teich(S)$. Of course, the identification $CC(N;S,X) \cong \Teich(S))$ also allows us to consider $V$ as a vector field on $CC(N;S,X)$. Conceptually this may be preferable as the hyperbolic structures determine $V$. However, much of the actual work (after the definition) will only involve Teichm\"uller space and we will move freely between the two viewpoints.
 As we will discuss below $V$ is the Weil-Petersson gradient of the negative of the {\em renormalized volume} function on $CC(N;S,X) \cong \Teich(S)$.

Our main result is the following:
\begin{theorem}
Let $(N;S)$ be relatively acylindrical and $M_t \in \Teich(S)$ be a flowline for $V = V_{(N;S,X)}$, then $M_t$ converges to $M_{\rm geod}$.
\label{theorem:main}
\end{theorem}

For the case of $N$ being acylindrical, the above states that the flow $V$ uniformizes $N$ in that every convex cocompact structure flows to the unique structure $M_{\rm geod}$ with totally geodesic boundary. This is also the structure with minimal convex core volume (see \cite{storm1}).

We note the existence of the manifold $M_{\rm geod}$ is a consequence of Thurston's hyperbolization theorem along with the deformation theory of Kleinian groups mentioned above and in fact the manifold $M_{\rm geod}$ only exists if $(N;S)$ is relatively acylindrical. More precisely the proof is (a special case of) the induction step in the proof of Thurston's theorem which is to find a fixed point of the {\em skinning map} on Teichm\"uller space. Thurston proved this by showing that the skinning map has bounded image (see \cite{Thurston:hype1}). McMullen gave an alternative proof by showing that this skinning map was a strict contraction (see \cite{McMullen:iter}). Our proof uses McMullen's contraction of the skinning map at two key moments  although for one, Thurston's bounded image theorem would also work. So we are not  giving a new proof of this existence theorem. It would be very interesting to give a more direct proof of Theorem \ref{theorem:main} that didn't depend on these two results which would give an alternative proof of the existence of $M_{\rm geod}$.

The deformation space $CC(N; S,X)$ is homeomorphic to an open ball and the vector field $V$ has a single, attracting zero (see \cite{moroianu:min,vargas:min}) so it may not seem surprising that the flow converges to this zero. However, the boundary of $CC(N; S, X)$ (appears to) exhibit fractal behavior which the vector find must wind its way through to find the zero. For example if $CC(N; S, X)$ is the a Bers slice $\bx$ then the Bers embedding identifies $\bx$ with a bounded  open topological ball  in the finite dimensional vector space $Q(X)$ of holomorphic quadratic differentials on $X$. When $\bx$ has complex dimension one Komori, Sugawa, Wada and Yamashita (see \cite{bers:images}) and Dumas (see \cite{dumas:bear} for images and \cite{dumas:bearprogram} for software) have drawn pictures that reveal this fractal behavior. More rigorously, also for dimension one,  Miyachi (see \cite{miyachi:cusp}) has shown that ``cusped'' manifolds on the boundary of $\bx$ correspond to cusps in the boundary of $\bx$ itself. By McMullen (\cite{mcmullen:cusps}) cusped manifolds are dense in the boundary of $\bx$ so together this implies that the boundary of $\bx$ has a dense set of cusps, a more concrete indication of the fractal nature of $\bx$. The flow $V$ is a natural flow that gives a contraction of these complicated domains to the Fuchsian basepoint.

The proof in the Bers slice case and the general case  differ only in that the general case requires  additional analysis to show that the extra components (called leopard spots) in the domain of discontinuity do not contribute to the limiting model flow detailed below. For clarity of exposition, we have isolated this additional analysis to Section \ref{spots}. 

We conclude this introduction with a informal discussion of the flow $V$ when $N$ is acylindrical (so $S = \partial N$). A construction of C. Epstein (see \cite{Epstein:surface}) describes a surface $Y'$ in $M_Y$ associated to the hyperbolic metric on $Y$. This surface cuts off a compact core of $M_Y$ which is closely related to the convex core. When the $L^2$-norm of the Schwarzian is small the curvature is small. In particular the integral of the mean curvature will be small. The difficulty is that this does not imply that the curvature is small everywhere in $Y'$ but only on the $\epsilon$-thick part of $Y'$ for some small $\epsilon$.  When we start flowing along $V$ in $CC(N)$ the flow will try to locally deform $Y'$ to decrease its curvature.
This puts the thick and thin parts in competition - to decrease the curvature in the thick part we need to increase it in the thin part and vice versa. Our central conclusion is that the thin part eventually wins - the flow will eventually decrease curvature in the thin part even at the cost of increasing it in the thick part. As this happen the short curves will become longer as the flow travels to a different point in Teichm\"uller space. This  process can repeat but  in \cite{wp-paper} we saw that it can only happen a finite number of times and eventually the Epstein surface will converge to a totally geodesic surface that bounds the convex core of the limiting manifold  which therefore must be $M_{\rm geod}$.

{\bf Acknowledgement.} We would like to thank MSRI  for their hospitality while portions of this work were being completed. We also thank David Dumas and Curt McMullen for helpful conversations. Finally we would like to especially thank the reviewer whose comments and suggestions greatly improved the paper.
 
 \subsection{The limiting model flow}\label{toy}
In order to prove our main theorem, we show that if a flowline does not converge to $M_{\rm geod}$ then we can extract a limiting model flow as the flowline tends a point in the Weil-Petersson completion. We then use the properties of this  model flow to obtain a contradiction. We now describe the limiting model flow.  

\begin{figure}[htbp] 
   \centering
   \includegraphics[width=3in]{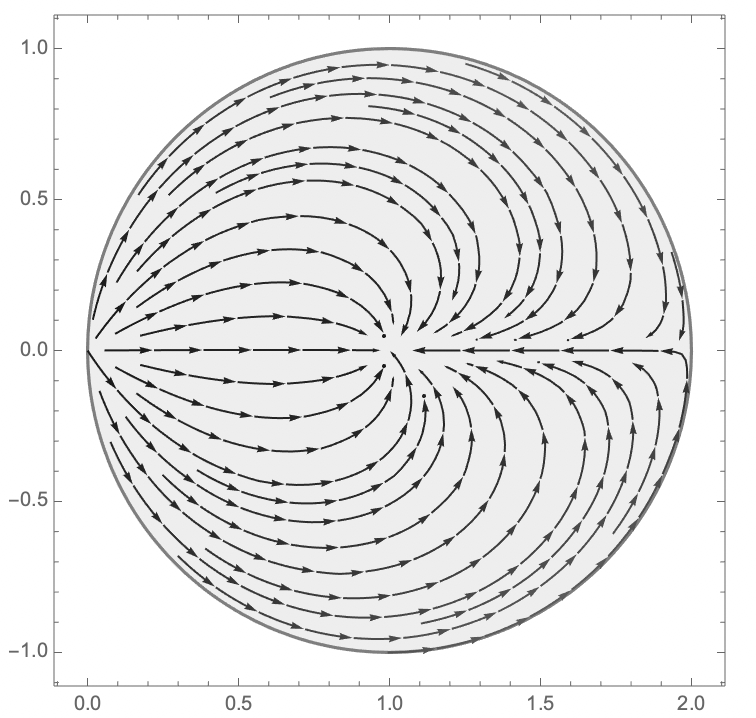} 
   \caption{Vector field $v$  on $\mathcal D$}
   \label{vectorfield}
\end{figure}

 Identifying the hyperbolic plane with the upper half-plane in $\CC$ we consider the space of univalent maps $f_c\colon \Hp \rightarrow \hat\CC$ of the form $f_c(z) = z^c$  whose image  is a Jordan domain. 
This is subspace of the space of all univalent maps on $\Hp$ and it corresponds to the open disk $\mathcal D = \{ f_c\ | \ |c-1| < 1\}$. (See Lemma \ref{univalent_c}.)
For each $f_c \in \mathcal D$ define a Beltrami differential $\mu_c$ on $\chat$ such
$$f^*_c \mu_c = -\frac{\overline{S(f_c)}}{\rho_{\Hp}}$$
and $\mu_c$ is zero on the complement of the image of $f_c$.

This family of Beltrami differentials defines a flow on $\cD$ as follows. For any $f_c \in \cD$ there is a family of quasiconformal homomeomorphisms  $\psi_t \colon \chat \to\chat$ whose infinitesimal Beltrami differential at $t=0$ is $\mu_c$ and there is a smooth path $f_{c_t}$ such that $f_{c_t}(\Hp) = \psi_t(f_c(\Hp))$. Furthermore the $\psi_t$ are defined for $t\in \R$ and the flowline starting at $f_{c_s} \in \cD$ is given by the formula $\psi_{t+s}\circ \psi_{-s}$.

The map $f_c\mapsto c$ is a homeomorphism form $\cD$ to the disk $|z-1| < 1$ in $\CC$ and the flow on $\cD$ induces the flow $c_t$ on the disk. In fact we have the following formula for the flow:
$$v(c) =\frac{1}{4}\left( |c|^4-2c\Re\left(c^2\right)-c^2+2c\right).$$
Although these last two paragraphs could be made rigorous as they are not necessary to prove our main result we will not do so. 
This formula for $v$ will follow from our derivation of the limiting vector field in our proof. However, to informally see the connection between this flow and the gradient flow $V$ one should view the domains $f_c(\Hp) \subset \chat$ as approximations for a component of the domain discontinuity of a Kleinian group where the imaginary axis is the axis of some short curve in the conformal boundary and the constant $c$ is the ratio of the complex length of the curve in the hyperbolic 3-manifold to its length on the boundary. There is an explicit formula for the derivative of this length ratio (see Section \ref{derivative_formula}) and in Theorem  \ref{lemma:c'limit} we calculate the limit of this formula as both the length on the boundary and the norm of the Schwarzian approach zero. The limiting formula we find is exactly the formula for the vector field $v$ above.

We can also connect the discussion here to our discussion of the flow in the previous section. In our model situation when $c$ is  near two corresponds to an acylindrical manifold where the surface $Y'$ has small curvature outside of the $\epsilon$-thin part where the closer $c$ is to two the smaller we can choose $\epsilon$. When $c$ is near one the corresponds to $Y'$ having small curvature everywhere. While the flow $v$ may initially appear to be converging to two it will eventually turn around and head towards one which corresponds $Y'$ being totally geodesic and the hyperbolic manifold $M_{\rm geod}$.

While our discussion here is only informal we will see that the properties of this limiting model flow play a crucial role in the proof of our main theorem.

\section{Weil-Petersson geometry}
Let $S$ be a closed surface of genus $g \geq 2$. Then Teichm\"uller space  $\Teich(S)$ is the 
space of marked conformal structures on $S$. Given $X \in \Teich(S)$ the cotangent space $T^*_X(\Teich(S))$  is $Q(X)$ the space of holomorphic quadratic differentials on $X$. We let $B(X)$ be the the space of Beltrami differentials on $X$. Then there is a pairing between $Q(X)$ and $B(X)$ given by
$$(\phi,\mu) = \int_X \phi\mu.$$
If we let $N(X) \subseteq B(X)$ be the annihilator of $Q(X)$ under this pairing, we obtain the identification $T_X(\Teich(S)) = B(X)/N(X)$.

  Given $\phi \in Q(X)$ and $z \in X$ then we define the pointwise norm by
$$\|\phi(z)\| = \frac{|\phi(z)|}{\rho_X(z)}$$
where $\rho_X$ is the hyperbolic metric on $X$. We define the $L^p$ norm of $\phi$, denoted $\|\phi\|_p$, to be the $L^p$ norm of the function $\|\phi(z)\|$ with respect to the hyperbolic area form on  $X$.
These $L^p$ norms define Finsler cometrics on the cotangent bundle of $\Teich(S)$ and dual Finsler metrics on the tangent bundle of $\Teich(S)$. When $p=2$ this norm comes from an inner product and therefore determines a Riemannian metric on $\Teich(S)$ called the  {\em Weil-Petersson metric}.  Classical results are that the Weil-Petersson metric is incomplete  (see \cite{Chu:noncompleteness,Wolpert:noncompleteness}) and  strictly negatively curved (see \cite{Tromba:sectional,Wolpert:sectional}).

In order to describe the Weil-Petersson completion $\overline{\Teich(S)}$, we first describe the {\em augmented Teichm\"uller space}. For further details on augmented Teichm\"uller space see \cite{Abikoff:degenerating},\cite{abikoff1980real} and \cite{Harvey:compactness}.

For $S$ a compact surface we let $\calC(S)$ be  the {\em complex of curves},  the simplicial complex organizing the isotopy classes of
simple closed curves on $S$ that do not represent boundary components.
To each isotopy class $\alpha$ we associate a vertex $v_\alpha$, and
each $k$-simplex $\sigma$ is the span of $k+1$ vertices whose
associated isotopy classes can be realized disjointly on $S$.

A point in the augmented Teichm\"uller space is given by a choice
of  multicurve $\tau$, a ($0$-skeleton of a) simplex in $\calC(S)$, and finite area hyperbolic structures on the
complementary subsurfaces $S \setminus \tau$. The elements of $\tau$ are the {\em nodes} and the point of the completion is a {\em noded Riemann surface}. Augmented Teichm\"uller space is
stratified by the simplices of $\calC(S)$: the collection of noded
Riemann surfaces with nodes determined by a given simplex $\sigma$
lies in a product of lower-dimensional Teichm\"uller spaces determined
by varying the structures on $S \setminus \tau$. This {\em
  stratum}, $\calS_\tau$, inherits a natural metric
from the Weil-Petersson metric, which by Masur (see \cite{Masur:WP})  is isometric to
the product of Weil-Petersson metrics on the Teichm\"uller spaces of
the complementary subsurfaces.

It follows by Masur also  that the augmented Teichm\"uller space is the Weil-Petersson completion $\overline{\Teich(S)}$. The completion naturally descends under the action of the
mapping class group to a finite diameter metric on the Deligne-Mumford
compactification of the moduli space of Riemann surfaces. The strata  of the completion can be described as follows; 
$$\calS_\tau = \{ X \in \overline{\Teich(S)}\ |\  \ell_\alpha(X) =0 \mbox{ if and only if } \alpha \in \tau\}$$
where $\ell_\alpha$ is the extended length function of $\alpha$.

\subsection{Length functions and Gardiner's Formula}\label{derivative_formula}
For an essential closed curve $\alpha$ on $S$ there are two natural length functions on the deformation space $CC(N;S,X) \cong \Teich(S)$. The first is just the usual length function $\ell_\alpha\colon \Teich(S) \to \R_+$ on Teichm\"uller space where $\ell_\alpha(Y)$ is the length of the geodesic representative of $\alpha$ on the hyperbolic surface $Y$. We have used this function already.

The {\em Gardiner formula} is a formula for the differential $d\ell_\alpha$. To state it we identify the universal cover of $Y$ with the upper half plane normalized so that the imaginary axis is an axis for $\alpha$. Then $z\mapsto e^{\ell_\alpha(Y)}z$ is an element of the deck group for $Y$. We let $A_\alpha$ be the quotient annulus for the action of this element. The annulus $A_\alpha$ also covers $Y$ so if $\mu$ is a Beltrami differential on $Y$ (representing a tangent vector in $T_Y\Teich(S)$) then $\mu$ lifts to a Beltrami differential $\mu_A$ on $A_\alpha$. The holomorphic quadratic differential $dz^2/z^2$ on $\Hp$ is invariant under the action of $z\mapsto e^{\ell_\alpha(Y)}z$ so descends to a Beltrami differential on $A_\alpha$. We will continue to refer to this quadratic differential as $dz^2/z^2$ on $A_\alpha$. If we let $\langle , \rangle_{A_\alpha}$ be the pairing between Beltrami differentials and quadratic differentials on $A_\alpha$ we have
\begin{theorem}[{Gardiner, \cite{Gardiner:theta}}]\label{gardiner}
The derivative of $\ell_\alpha$ on $\Teich(S)$ is given by the formula
$$d\ell_\alpha(\mu) = \frac2\pi\left\langle \mu_A, \frac{dz^2}{z^2}\right\rangle_{A_\alpha}.$$
\end{theorem}

If we consider $\alpha$ as a closed curve in the hyperbolic 3-manifold $M_Y$ then $\alpha$ has a {\em complex length}. The real part is just the length of the geodesic representative of $\alpha$ in $M_Y$ while the imaginary part measures the twisting along the geodesic. We need for the imaginary part to be a well defined real number (rather than just a number mod $2\pi$). For this reason the definition is somewhat involved.

Let $\Gamma_Y$ be a Kleinian group uniformizing $M_Y \in CC(N;S,X)$ and let $\Omega_Y$ be the components of the domain of discontinuity of $\Gamma_Y$ that cover $Y$. The pre-image of $\alpha$ in $\Omega_Y$ will be a collection of arcs. Fixing an orientation of $\alpha$ fixes an orientation of each of the arcs and we can assume that the pre-image contains an oriented arc $\tilde\alpha$ with initial endpoint $0$ and terminal endpoint $\infty$. Let $\Omega\subset\Omega_Y$ be the component of the domain of discontinuity that contains $\tilde\alpha$ and let $\log$ be a branch of the logarithm defined on $\Omega$. We then define
$$\cL_\alpha\colon CC(N;S,X) \to \CC$$
by
$$\mathcal L_\alpha(M_Y) = \log \tilde\alpha(\ell_\alpha(Y)) - \log \tilde\alpha(0)$$
where we are assuming the that $\tilde\alpha$ is parameterized by arc length. The complex length $\cL_\alpha(M_Y)$ is independent of the choices we made in its definition.

The complex length is a holomorphic function on $CC(N:S,X)$ and we would also like a formula for its derivative. Miyachi observed (\cite[First proposition, Section 8]{miyachi:cusp}) that the proof of Theorem \ref{gardiner} can also be applied to find the derivative of the complex length.  As $CC(N;X,S) \cong \Teich(S)$ if $\mu$ is a tangent vector in $T_Y\Teich(S)$ we can also consider it as a tangent vector to $CC(N;X,S)$ at $M_Y$. Then $\mu$ will lift to a $\Gamma_Y$-invariant Beltrami differential $\tilde\mu$ on $\Omega_Y$. We can extend $\tilde\mu$ to be zero everywhere else. The Kleinian group $\Gamma_Y$ will contain the element $z\mapsto e^{\cL_\alpha(M_Y)}z$. The quotient of the $\CC\smallsetminus\{0\}$ under the action of this element will be a torus $T_\alpha$ and the Beltrami differential $\tilde\mu$ will descend to a Beltrami differential $\mu_T$ on $T_\alpha$. The quadratic differential $dz^2/z^2$ will also descend to a quadratic differential on $T_\alpha$. If $\langle, \rangle_{T_\alpha}$ is the pairing on $T_\alpha$ we have
\begin{theorem}[{Miyachi, \cite{miyachi:cusp}}]\label{miyachi}
The differential of $\cL_\alpha$ is given by the formula
$$d\cL_\alpha(\mu) = \frac1\pi\left\langle \mu_T, \frac{dz^2}{z^2}\right\rangle_{T_\alpha}.$$
\end{theorem}

\begin{figure}[htbp] 
   \centering
   \includegraphics[width=1.5in]{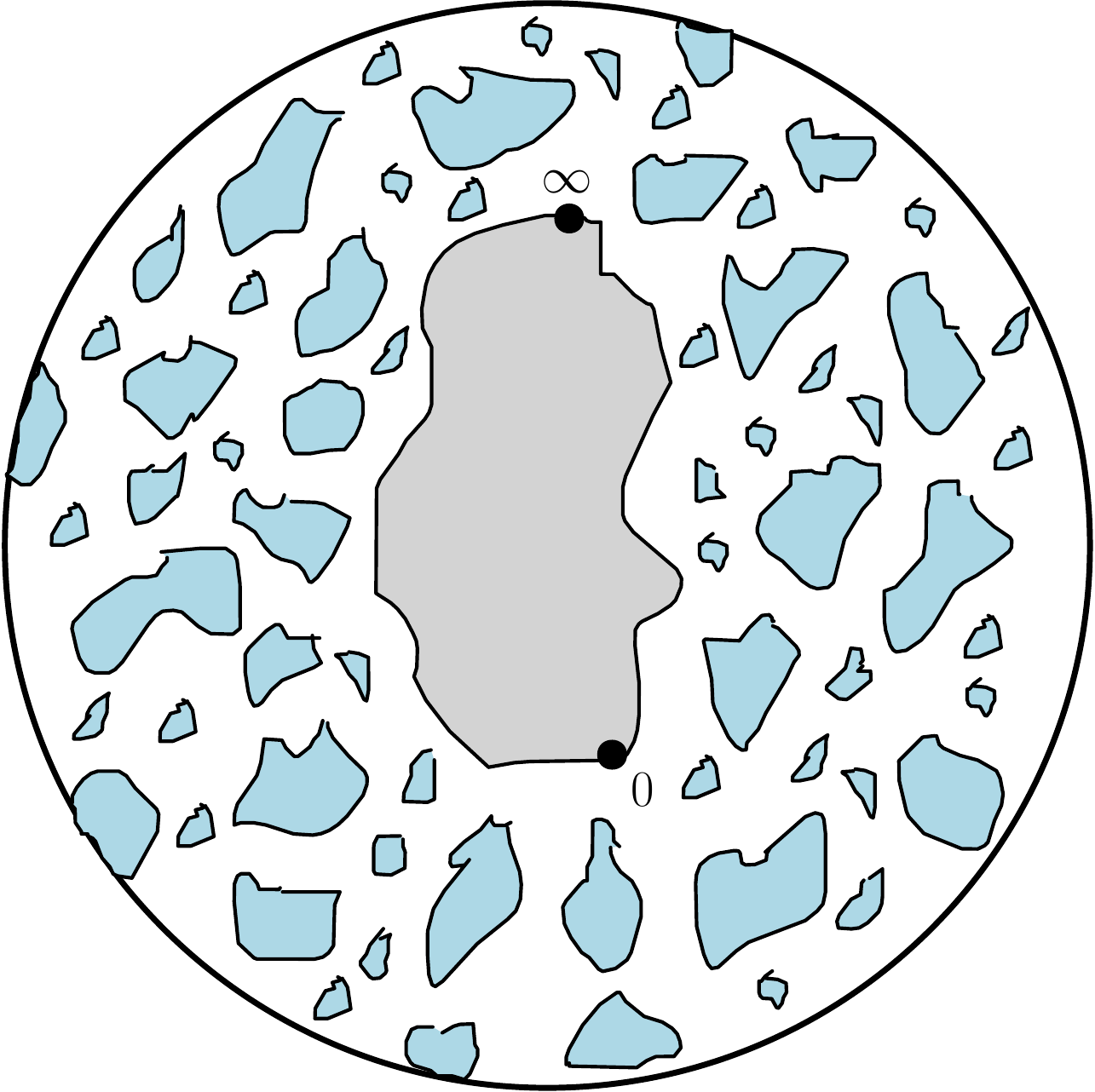}    \hspace{.5in}\includegraphics[width=1.5in]{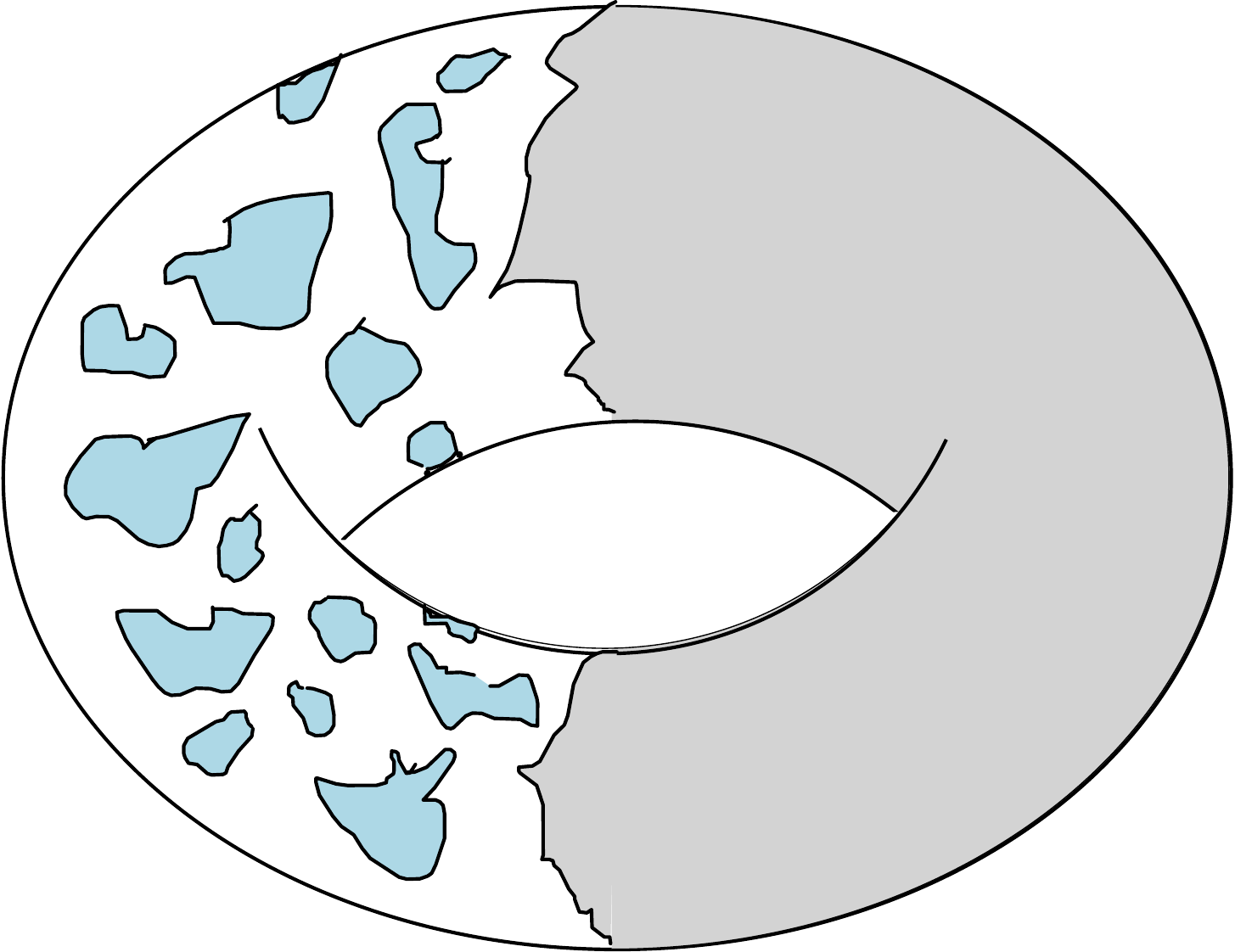} 
   \caption{Leopard spots on Riemann sphere and quotient torus}
   \label{leopard}
\end{figure}

For what we will do below it will be useful to decompose this second Gardiner formula into two distinct terms. To describe this decomposition we note that the image of $\Omega$ in $T_\alpha$ will be an essential annulus while every other component of $\Omega_Y$ will map homeomorphically into $T_\alpha$ (see figure \ref{leopard}). We then write
$$\mu_T = \mu_T^{\rm cen} + \mu_T^{\rm aux}$$
where $\mu_T^{\rm cen}$ has support on the image of $\Omega$ and the support of $\mu_T^{\rm aux}$ is in the image of the other components of $\Omega_Y$. Then
$$\frac1\pi\left\langle \mu^{\rm cen}_T, \frac{dz^2}{z^2}\right\rangle_{T_\alpha} \qquad \mbox{and} \qquad \frac1\pi\left\langle \mu^{\rm aux}_T, \frac{dz^2}{z^2}\right\rangle_{T_\alpha} $$
are the {\em central} and {\em auxillary} terms of the differential $d\cL_\alpha(\mu)$.

It will also be useful to write the central term as a pairing on the annulus $A_\alpha$. For this let
$$g\colon \Hp\to \Omega$$
be the uniformizing map, normalized so that $g$ takes the imaginary axis to $\tilde\alpha$. Then the pull backs $g^*\tilde\mu$ and $g^*(dz^2/z^2)$ are both invariant under the isometry $z\mapsto e^{\ell_\alpha(Y)}z$ and descend to objects on $A_\alpha$.  In fact the Beltrami differential will be the Beltrami differential $\mu_A$ that we defined above. We can then write the central term as
$$\frac1\pi\left\langle \mu^{\rm cen}_T, \frac{dz^2}{z^2}\right\rangle_{T_\alpha} = \frac1\pi\left\langle \mu_A, g^*\left(\frac{dz^2}{z^2}\right)\right\rangle_{A_\alpha}.$$

Note that both the central and auxillary terms only depend on the Beltrami differential $\mu$. For later convenience we let $J_\alpha(\mu)$ be the auxillary term. Note that for  a Bers slice  $\Omega$ is the only component of $\Omega_Y$ so the auxillary term is always zero. The extra work in the relatively acyndrical case is estimating $J_\alpha(\mu)$ when $\mu$ is a harmonic Beltrami differential with small $L^2$-norm.

\section{Limits of the flow}
As mentioned in the introduction, the flow $V$ is the negative of the Weil-Petersson gradient of the {\em renormalized volume} function ${\Rvol}$ on $CC(N;S,X)$. The definition of renormalized is somewhat involved and not necessary for our work here. Therefore we will omit it and restrict ourselves to discussing some of its important properties. The first of these is a varitional formula:

\begin{theorem}[\cite{ZT,TT,KS12}]
\label{variational}
For tangent vector $\mu \in   T_{M_Y} CC(N) \cong T_Y\Teich(S)$ we have
$$d\Rvol(\mu) = \Re \int_{Y} \phi_{Y}\mu.$$
\end{theorem}
This formula implies that our flow $V$ is the Weil-Petersson gradient of $-{\Rvol}$.

By the Nehari bound on the norm of the Schwarzian derivative of a univalent map (see \cite{Nehari:schwarzian}), $V$ is bounded with respect to the Teichm\"uller metric on $\Teich(S)$. Therefore as the Teichm\"uller metric is complete, the flowlines exists for all time (see \cite{BBB} for further details).

Also by the gradient description of $V$ it follows that along a flowline a $M_t$
$${\Rvol}(M_0)-{\Rvol}(M_T) = \int_{0}^T \|\phi_{Y_t}\|_2^2 dt.$$
As ${\Rvol} \geq 0$ (see \cite{BBB}), it follows that 
$$\int_0^\infty  \|\phi_{Y_t}\|_2^2 dt < \infty.$$

  We now describe the further properties of the flow proved in \cite{wp-paper} that we will need in our analysis.

\begin{theorem}[{Bridgeman-Brock-Bromberg, \cite{wp-paper}}]
Let  $(N,S)$ be relatively acylidrical and $M_t = (X,Y_t)$ be a  flowline for $V$ on $CC(N;S,X)$. Then
\begin{enumerate}
\item $Y_t \rightarrow \hat Y \in \overline{\Teich(S)}$ in the Weil-Petersson completion. Thus $\hat Y$ is a noded Riemann surface.
\item $\|\phi_{Y_t}\|_2 \rightarrow 0$ as $t\rightarrow \infty$.
\end{enumerate}
\label{flowfacts}
\end{theorem}
In order to prove our main theorem, we need to prove that  the set of nodes of any limit $\hat Y$ is empty or alternately that $\hat Y \in \Teich(S)$. We will do this by assuming $\hat Y$ is noded and consider the limits of the projective structure as we zoom in on the nodes. This will produce our limiting model flow which will allow us to obtain a contradiction.

\section{Taking limits at a node}
The  length functions $\ell_\alpha$  and $\cL_\alpha$ are smooth functions on $CC(N;S,X) \cong \Teich(S)$ if we pair their differentials against the vector field $V$ we get a function on $\Teich(S)$. We would like to take the limit of these functions along sequence $Y_n$ in $\Teich(S)$ where $\|V(Y_n)\|_2 \to 0$. To do this we will embed $CC(N;S,X)$ in a space of univalent functions. We then use normal families theorems for holomorphic functions to evaluate the limit. We begin by setting up our space.

Let 
$$\cU =\{ \phi | \mbox{ $\phi= S(f)$ for some univalent $f\colon \Hp \to \CC$}\}$$
 be the space of quadratic differentials that arise as Schwarzian derivatives of univalent functions from $\Hp$ (realized as the upper half plane) to $\CC$. We give $\cU$ the compact-open topology and recall some basic facts:
 \begin{itemize}
 \item $\cU$ is compact.
 
 \item If $\phi_n\to \phi$ and $f_n\to f$ (in the compact-open topology) with $S(f_n) = \phi_n$ then either $f$ is constant or $f$ is univalent and $S(f) = \phi$.
 \end{itemize}
 
 The following result is essentially Montel's theorem.
\begin{theorem}\label{normal}
Let $z_0$, $z_1$ and $z_2$ be distinct points in $\overline{\mathbb H}^2$ and $w_0$, $w_1$ and $w_2$ distinct points in $\chat$. Let $\cF$ be a family of holomorphic maps on $\Hp$ that extend continuously to the $z_i$ if they are on the boundary of $\Hp$ and assume that for all $f\in\cF$ we have $f(z_i) = w_i$ for $i=0,1$ and $2$ while if $z$ is not in $\{z_0, z_1, z_2\}$ then $f(z)$ is not in $\{w_0, w_1,w_2\}$. Then $\cF$ is a normal family. If $f$ is a limit of a sequence in $\cF$  then $f(z_i) = w_i$ for $z_i \in \Hp$ and if $f$ is non-constant then  $f(z_i) = w_i$ for all $i=0,1,2$.
\end{theorem}

For  $s \in \R$ we say that $\phi\in \cU$ is {\em $s$-invariant} if it is invariant under the isometry $z\mapsto e^sz$ (as a quadratic differential). Let $\cU_\Z\subset \cU$ be the subspace of quadratic differentials that are $s$-invariant for some $s \neq 0$. We can define a function 
$$\ell\colon \cU_\Z\to [0,\infty)$$
by taking $\ell(\phi)$ to be the infimum over all positive $s$ where $\phi$ is $s$-invariant. Note that if $s_n \to s$ and $\phi_n$ are $s_n$-invariant then if $\phi_n\to \phi$ we have that $\phi$ is $s$-invariant.  This implies that $\ell$ is continuous on $\cU_\Z$. It is possible that $\ell(\phi) =0$ and we let $\cU^0_\Z = \ell^{-1}(0)$. This space will be of particular interest.

\begin{lemma}\label{univalent_c}
If $\phi \in \cU^0_\Z$ then 
$$\phi(z) = \left(\frac{1-c^2}{2}\right)\frac{dz^2}{z^2}$$
 with $|c-1|\le 1$. If $c\neq 0$ let $g_c(z) = z^c/i^c$. If $c= 0$ let $g_c = \log z$. Then $\phi = S(g_c)$.
\end{lemma}

{\bf Proof:} If $\ell(\phi) = 0$ then $\phi$ is $s_n$-invariant for a sequence of $s_n>0$ with $s_n\to 0$. If $\phi$ is $s$-invariant then $\phi$ is $sk$-invariant for all $k \in \Z$. Together this implies that $\phi$ is $s$-invariant for a dense set of $s \in \R$ so by continuity $\phi$ is $s$-invariant for all $s\in\R$. From this invariance we see that if $\phi(i) = -C dz^2$ for some $C \in \CC$ then $\phi(it) = -\frac C{t^2} dz^2$ for all $t>0$. As $\phi$ is holomorphic this implies that $\phi(z) = \frac{C}{z^2} dz^2$.

A direct computation gives that $S(g_c) = \frac{1-c^2}{2}\cdot\frac{dz^2}{z^2}$. We'll show that $g_c$ is univalent exactly when $|c-1|\le 1$. We can write $g_c(z) = e^{c\log z}/i^c$. Let $\cS$ be the open horizontal strip in $\CC$ between the lines $\Im z = 0$ and $\Im z = \pi$ and let $\cS_c = \{cz\ |\ z\in \cS\}$. Then the image of  $\Hp$ (as the upper half plane) under the map $z\mapsto c\log z$ is $\cS_c$. The exponential map restricted to $\cS_c$ is injective exactly when vertical lines in $\CC$ intersect $\cS_c$ in a segment of length $< 2\pi$. The intersection of a vertical line with $\cS_c$ has length $|c|\pi/\cos(\theta)$ where $\theta = \arg(c)$. Thus we have univalence if  $|c| \leq 2\cos(\theta) = 2\Re(c)/|c|$ giving $|c|^2-2\Re(c) \leq 0$. Completing the square, we get $|c-1|^2\leq 1$.
\eproof

The invariance of the quadratic differentials $\phi\in\cU_\Z$ imply invariance for any univalent map $f$ with $S(f) = \phi$. The next lemma makes this precise.
\begin{lemma}\label{ell_0}
If the quadratic differential $\phi\in \cU_\Z$ is $s$-invariant and $f\colon \Hp\to \CC$ is univalent with $S(f) = \phi$ then there exists a $\psi \in \psl$ with $\psi\circ f(z) = f(e^s z)$. Furthermore $\psi$ is either loxodromic or parabolic and $f$ extends  continuously to $0$ and $\infty$ with $f(0)$ the repelling fixed point of $\psi$ and $f(\infty)$ the attracting fixed point if $\psi$ is loxodromic and $f(0) = f(\infty)$ if $\psi$ is parabolic.
\end{lemma}

{\bf Proof:} If $\phi$ is $s$-invariant and $f$ is a univalent map with $\phi= S(f)$ then the map $z\mapsto f(e^s z)$ also has  Schwarzian $\phi$. As two maps with the same Schwarzian differ by post-composition of an element of $\psl$ we have that there exists a $\psi \in \psl$ with $\psi \circ f(z) = f(e^s z)$. Iterating this formula, for any positive integer $k$ we have that $\psi^k\circ f(z) = f(e^{sk}z)$. This implies that $f(\Hp)$ is $\psi$-invariant and that the action of $\psi$ on the simply connected space $f(\Hp)$ doesn't have fixed points so $\psi$ must be loxodromic or parabolic. Furthermore the equation $\psi^k\circ f(z) = f(e^{sk}z)$ implies that $f$ extends continuously to $0$ and $\infty$ with $f(0)$ the repelling fixed point of $\psi$ and $f(\infty)$ the attracting fixed point. \eproof

The element $\psi$ is the {\em holonomy} of $\phi$. Let $\cU^+_\Z \subset \cU_\Z$ be the subspace of quadratic differentials $\phi$ where $\ell(\phi) >0$ and the holonomy is loxodromic. For $\phi \in \cU^+_\Z$ we define a complex length function $\cL\colon \cU^+_\Z\to \CC$ as follows. Let $f\colon \Hp\to \CC$ be a univalent map with $f(0) = 0$, $f(\infty) = \infty$ and $S(f) = \phi$. Choose a logarithm, $\log_\phi$, on $f(\Hp)$. We then define $\cL(\phi) = \log_\phi\left(f(e^{\ell(\phi)}z)\right) - \log_\phi(f(z))$. The expression on the right is independent of the choice of $f$, the choice of $z \in \Hp$ and the choice of logarithm.

We also define a function $c\colon \cU_\Z^+\to \CC$
by 
$$c(\phi) = \frac{\cL(\phi)}{\ell(\phi)}.$$

\begin{lemma}\label{continuous_c}
The function $c$ extends continuously to $\cU^0_\Z$ and for $\phi \in \cU^0_\Z$ we have $$\phi(z) = \frac{1-c(\phi)^2}{2}\cdot \frac{dz^2}{z^2}.$$
\end{lemma}

{\bf Proof:} If $f$ is univalent and $S(f)=\phi \in \cU^+_\Z$ then we have the Bers' inequality (or McMullen's interpretation of the Bers' inequality, \cite[Proposition 6.4]{McMullen:iter}):
\begin{equation}\label{bersineq}
\frac{1}{\ell(\phi)} \leq \frac{2\Re \mathcal L(\phi)}{|\mathcal L(\phi)|^2}.
\end{equation}
The statement is usually made in the context of quasifuchsian groups but the proof goes through without change in our setting. From this we see
\begin{eqnarray*}
|c(\phi)|^2 = \frac{|\mathcal L(\phi)|^2}{\ell(\phi)^2} &\leq& 2\Re\left(\frac{\mathcal L(\phi)}{\ell(\phi)}\right) = 2\Re(c(\phi)),\\
|c(\phi)|^2-2\Re(c(\phi)) &\leq& 0, \\
 |c(\phi) - 1|^2-1 &\leq& 0.
\end{eqnarray*}
Thus $|c(\phi)-1| \leq 1$.

Now assume $\phi \in \cU_\Z^0$. Then by Lemma \ref{ell_0} we have $\phi = C\frac{dz^2}{z^2}$ for some $C \in \CC$. Take a sequence $\phi_n$ in $\cU_\Z^+$ with $\phi_n \to \phi$ and let $\ell_n = \ell(\phi_n), \cL_n = \cL(\phi_n)$ and $c_n = c(\phi_n)$. Since $\ell$ is continuous on $\cU_\Z$ we have that $\ell_n\to 0$ and we can choose integers $k_n$ such that $k_n\ell_n\to 1$. We also fix univalent maps $f_n$ with $S(f_n) = \phi_n$ and the normalization $f_n(\infty) = \infty$, $f_n(i) =0$ and $f_n(ie) = 1$.  By Theorem \ref{normal} the $f_n$ will be a normal family. Let $f$ be a limit of a subsequence of the $f_n$. Again by Theorem \ref{normal}, $f(i) =0$ and $f(ie) =1$ so the limit won't be constant  which in turn implies that $f(\infty) = \infty$ and that $f$ is univalent. It follows that $S(f_n) = \phi_n \to \phi = S(f)$. As $f$ is the only map with $S(f) = \phi$ plus the given normalization, any convergent subsequence of the $f_n$ will converge to $f$. Therefore $f_n$ converges to $f$ (uniformly on compact sets) without passing to a subsequence.

If $\psi_n$ is the element of $\psl$ with $f_n(e^{\ell_n}z) = \psi_n\circ f_n(z)$ then the attracting fixed point of $\psi_n$ and all its powers is $\infty$. Therefore $\psi_n^{k_n}(z)= e^{k_n \cL_n} z + A_n$ for some $A_n \in \CC$. As $c_n$ lies in a compact set we can pass to a subsequence such that $c_n \to b$ and therefore $k_n\cL_n = k_n\ell_n c_n \to b$.  Taking the limit of the equation $f_n(e^{k_n\ell_n}z) = e^{k_n\cL_n}f_n(z) + A_n$ we have
$$f(ez) = e^b f(z) + \lim_{n\to\infty} A_n.$$
Substituting in $i$ for $z$ we see that $A_n\to 1$. It follows that $e^b$ is determined by $f$ (and hence $\phi$). This determines $b$ up to a multiple of $2\pi i$. However, we also have $|b-1| \le 1$ which implies that $-1\le \Im b \le 1$ so the equation uniquely determines $b$. In particular, $c_n\to b$ before passing to a subsequence.

Let $g_c$ be the the univalent map given by Lemma \ref{univalent_c} with $S(g_c) = \phi$. Then $g_c = \beta\circ f$ for some $\beta\in\psl$. 
Note that $\psi^{k_n}_n \to \psi$ in $\psl$ with $\psi(z) = e^bz + 1$.  If $b=0$ then $f(0) = f(\infty) = \infty$ so we must have that $g_c(0) = g_c(\infty)$ and this only occurs when $c=0$ and $g_c(z) = \log z$. If $b\neq 0$ then the repelling and attracting fixed points of $\beta\circ\psi\circ\beta^{-1}$ are $g_c(0) = 0$ and $g_c(\infty) = \infty$. This implies that $\beta\circ\psi\circ\beta^{-1}(z) = e^b z$ and $g_c(ez) = e^bg_c(z)$. It follows that $b=c$. In both cases we can define $c(\phi) = \underset{n\to\infty}\lim c_n$ and we have that $c$ extends continuously to $\cU_\Z^0$ with
$$\phi(z) = \frac{1-c(\phi)^2}2\cdot \frac{dz^2}{z^2}.$$
\eproof

We now describe the map from $CC(N;S,X)$ to $\cU^+_\Z \subset \cU$. It will depend on a choice of essential closed curve $\alpha$ on $S$. Given $M_Y \in CC(N;S,X)$ let $\Omega$ be a component of the domain of discontinuity that projects to the component of $Y$ that contains $\alpha$. Let $f_Y\colon \Hp\to \Omega$ be a uniformizing univalent map that takes the imaginary axis in $\Hp$ to an axis for $\alpha$ and let $\phi_Y = S(f_Y)$. While we have made several choices the quadratic differential $\phi_Y$ is independent of our choices. We can then define $\Psi_\alpha\colon CC(N;S,X) \to \cU$ by $\Psi_\alpha(M_Y) = \phi_Y$. This map is continuous and $\ell_\alpha = \ell\circ\Psi_\alpha$. Since $\ell_\alpha$ is positive on $CC(N;S,X)$ we have that the image of $\Psi_\alpha$ lies in $\cU^+_\Z$.

Given $\phi \in \cU$ we define 
$$\mu_\phi = -\frac{\bar\phi}{\rho_\Hp}$$ 
Thus $\mu_\phi$ is the negative of the corresponding harmonic Beltrami differential. We also let $A_s$ be the quotient of $\Hp$ under the action $z\mapsto e^s z$ and $\langle, \rangle_s$ the pair of quadratic differentials and Beltrami differentials on $A_s$. If $\phi$ is $s$-invariant as a quadratic differential then $\mu_\phi$ is $s$-invariant as a Beltrami differential and descends to a Beltrami differential on $A_s$. We define a function
$$F_\ell\colon \cU_\Z^+ \to \R$$
by $$F_\ell(\phi) = \frac{2}{\pi\ell(\phi)} \Re \left\langle \mu_\phi, \frac{dz^2}{z^2}\right\rangle_{\ell(\phi)}.$$
As an immediate consequence of Gardiner's formula (Theorem \ref{gardiner}) we have:
\begin{lemma}\label{length_der}
$$d\log\ell_\alpha(V(M_Y)) = F_\ell(\Psi_\alpha(M_Y))$$
\end{lemma}

Therefore, to study the continuity of $d\ell_\alpha(V)$ we will study the continuity of $F_\ell$ on $\cU_\Z$.

\begin{lemma}\label{Fell_cont}
The function $F_\ell$ extends continuously to $\cU^0_\Z$ with
$$F_\ell(\phi) =\frac{\Re( c(\phi)^2)-1} 2$$
for $\phi \in \cU^0_\Z$.
\end{lemma}

{\bf Proof:} The key to the proof is that if $\mu$ and $\psi$ are $s$-invariant then $\langle \mu, \psi\rangle_{sk} = k\langle \mu, \psi\rangle_s$ for all positive integers $k$. Then if $\phi_n \to \phi$ for $\phi \in \cU^0_\Z$ we can choose $k_n$ such that $k_n\ell(\phi_n) \to 1$ so that
\begin{eqnarray*}
F_\ell(\phi_n) & = & \frac{2}{\pi\ell(\phi)} \Re\left\langle \mu_{\phi_n}, \frac{dz^2}{z^2} \right\rangle_{\ell(\phi_n)} \\
& = & \frac2\pi \Re\left\langle \frac{1}{k_n\ell(\phi_n)} \mu_{\phi_n}, \frac{dz^2}{z^2}\right\rangle_{k_n\ell(\phi_n)} \\
 &\longrightarrow & \frac2\pi \Re\left\langle \mu_\phi, \frac{dz^2}{z^2} \right\rangle_1.
\end{eqnarray*}
This shows the $F_\ell$ extends continuously to $\cU_\Z^0$.

For $\phi \in \cU^0_\Z$ we have $\phi(z) = \frac{1-c(\phi)^2}2 \cdot\frac{dz^2}{z^2}$. The pairing is easier to calculate in the strip model for $\Hp$ (the region $\cS$ between the horizontal lines $\Im z =0$ and $\Im z= \pi$) with area form $\rho_\Hp = 1/\sin^2 y$. In this model the quadratic differential $dz^2/z^2$ becomes $dz^2$ so $\mu_\phi(z) = \sin^2y (\overline{c(\phi)}^2-1)/2$ and
\begin{eqnarray*}
\left\langle \mu_\phi, \frac{dz^2}{z^2} \right\rangle_1 &=& \frac{\overline{c(\phi)^2}-1}2\int_0^\pi\int_0^1 \sin^2 y dx dy\\
&=& \frac{\pi\left(\overline{c(\phi)^2}-1\right)}4.
\end{eqnarray*}
Taking the real part and multiplying by $2/\pi$ gives the claimed formula for $F_\ell(\phi)$. \eproof

We would like to similarly define a function on $\cU^+_\Z$ for the differential $d\log\cL_\alpha$. We will not be able to do this exactly but instead give a formula for the central term. We will need to evaluate the auxillary term separately.

Given $\phi\in \cU^+_\Z$ let $g_\phi\colon \Hp\to \CC$ be the univalent map with $S(g_\phi) = \phi$ and $g_\phi(0) =0, g_\phi(\infty) =\infty$ and $g_\phi(i) = i$. By Theorem \ref{normal}, this is a normal family. We observe that if $\phi$ is $s$-invariant the quadratic differential $(g_\phi)^* \left(\frac{dz^2}{z^2}\right)$ is $s$-invariant. We then define a function
$$F_\cL\colon \cU_\Z^+\to \CC$$
by
$$F_\cL(\phi) = \frac1{\pi\cL(\phi)} \left\langle \mu_\phi, (g_\phi)^*\left(\frac{dz^2}{z^2}\right)\right\rangle_{\ell(\phi)}.$$
From Theorem \ref{miyachi} and the discussion following it we have:
\begin{lemma}\label{three_length_der}
$$d\log\cL_\alpha(V(M_Y)) = F_{\cL}(\Psi_\alpha(M_Y)) + \frac{J_\alpha(V(M_Y))}{\cL_\alpha(M_Y)}.$$
\end{lemma}

We now establish the continuity of $F_\cL$.
\begin{lemma}\label{FL_cont}
The function $F_\cL$ extends continuously to $\cU_\Z^0$ with 
$$F_\cL(\phi) = \frac{c(\phi)\left(\overline{c(\phi)}^2-1\right)}{4} .$$
\end{lemma}

{\bf Proof:} The proof is similar to Lemma \ref{Fell_cont}. Assume that $\phi_n$ in $\cU^+_\Z$ converges to $\phi\in\cU_\Z^0$ and let $g_n = g_{\phi_n}$, $\ell_n = \ell(\phi_n)$, etc. Choose integers $k_n$ with $k_n\ell_n \to 1$. By Lemma \ref{continuous_c} the function $c$ is continuous and therefore $k_n\cL_n=k_n\ell_n c_n \to c(\phi)$.

In Lemma \ref{gn_converge} below we'll show that
$$\frac1{c^2_n} g_n^* \left(\frac{dz^2}{z^2}\right) \to \frac{dz^2}{z^2}.$$
We assume this for now and then as in in Lemma \ref{Fell_cont} we have 
\begin{eqnarray*}F_{\cL}(\phi_n) & = & \frac1\pi\left\langle\frac1{k_n\cL_n}\mu_n, (g_n)^*\left(\frac{dz^2}{z^2}\right)\right\rangle_{k_n\ell_n} \\ & = & \frac{c_n}\pi\left\langle\frac1{k_n\ell_n}\mu_n, \frac{1}{c^2_n}(g_n)^*\left(\frac{dz^2}{z^2}\right)\right\rangle_{k_n\ell_n} \\ 
&\longrightarrow & \frac{c(\phi)}{\pi}\left\langle \mu_\phi, \frac{dz^2}{z^2}\right\rangle_1\\
&=& \frac{c(\phi)\left(\overline{c(\phi)^2}-1\right)}4.
\end{eqnarray*}
\eproof

\begin{lemma}\label{gn_converge}
$$\frac1{c^2_n} g_n^* \left(\frac{dz^2}{z^2}\right) \to \frac{dz^2}{z^2}$$
\end{lemma}

{\bf Proof:} We first assume that $c(\phi) \neq 0$. By our normalization $g_n(e^{k_n\ell_n} z) =e^{k_n\cL_n} g_n(z)$ so (after possibly passing to a subsequence) as $n\to \infty$ we have $g(ez) = e^{c(\phi)} g(z)$ where $g_n \to g$. Since $c(\phi)\neq 0$ we have that $g$ is non-constant and it follows that $g$ fixes $0$, $i$ and $\infty$ and $S(g) = \phi$. This implies that $g(z)= e^{c(\phi)}/i^{c(\phi)}$. It follows that
$$(g_n)^*\left(\frac{dz^2}{z^2}\right) \to g^*\left(\frac{dz^2}{z^2}\right) = c(\phi)^2\frac{dz^2}{z^2}$$ 
so the lemma holds in the case that $c(\phi) \neq 0$.

When $c(\phi) =0$ (and therefore by Lemma \ref{univalent_c} $\phi = \frac12\cdot\frac{dz^2}{z^2}$)   it will be necessary to choose a different normalization for the $g_n$ so that they don't converge to a constant function. In particular, similar to  the proof of Lemma \ref{continuous_c} we choose univalent functions $f_n$ with $f_n(\infty) = \infty$, $f_n(i) = i\pi$ and $f_n(ie) = 1+i\pi$ and whose Schwarzian is $\phi_n$. By Theorem \ref{normal}, the $f_n$ form a normal family. As  any limiting function will be non-constant the Schwarzians $\phi_n$ will also converge to the Schwarzian of the limit. As by assumption $\phi_n \to \phi$ we have that any limiting function has Schwarzian $\phi$. The chosen normalizations of $f_n$ will also persist in the limit. Together these conditions imply that the only possible limiting function is $f(z) = \log z$ so $f_n$ converges to $f(z) = \log z$ uniformly on compact sets.

Also, as in Lemma \ref{continuous_c}, we have the equation $f_n\left(e^{k_n\ell_n}z\right) = e^{k_n\cL_n} f_n(z) + A_n$ with $A_n \to 1$. If we choose $\beta_n \in \psl$ with $g_n = \beta_n \circ f_n$ then $g_n^*(dz^2/z^2) = f_n^* (\beta_n^*(dz^2/z^2))$. To calculate $\beta_n^*(dz^2/z^2)$ we observe that $\beta_n$ takes the attracting and repelling fixed points of $z\mapsto e^{k_n\cL_n} z + A_n$ to the attracting and repelling fixed points of $z\mapsto e^{k_n\cL_n} z$. In particular $\beta_n\left(A_n/\left(1-e^{k_n\cL_n}\right)\right) = 0$ and $\beta_n(\infty) = \infty$. It follows that
$$\beta_n^*\left(\frac{dz^2}{z^2}\right) = \frac{dz^2}{\left(z - \frac{A_n}{1-e^{k_n\cL_n}}\right)^2}.$$
As $k_n\cL_n = c_n k_n \ell_n \to c(\phi)\cdot 1 = 0$ this quadratic differential will converge to zero. However, if we multiply the denominator by $(k_n\cL_n)^2$ the denominator will converge to $1$ (since $A_n \to 1$) and the quadratic differential will converge to $dz^2$. We further have that $k_n\cL_n/c_n \to 1$ and it follows that if we
 divide by $c^2_n$ we have
$$\frac1{c^2_n} \beta_n^*\left(\frac{dz^2}{z^2}\right) \to dz^2$$
and since $f_n \to \log z$ this gives
$$\frac1{c_n^2} g_n^* \left(\frac{dz^2}{z^2}\right) =  f_n^*\left(\frac{1}{c_n^2} \beta_n^*\left(\frac{dz^2}{z^2}\right)\right)\to \frac{dz^2}{z^2}.$$
This proves the lemma when $c(\phi) = 0$.
\eproof

We can now prove our limiting formulas for the derivatives of $\ell_\alpha$, $\cL_\alpha$ and $c_\alpha$.
\begin{theorem}\label{lemma:c'limit}
Let $(N;S)$ be a relatively acylindrical pair and $\alpha$ an essential simple closed curve in $S$. Let $Y_n$ be a sequence in $\Teich(S)$ with the volume of $C(M_{Y_n})$ bounded, $\|V(Y_n)\|_2 \to 0$, $\ell_\alpha(Y_n) \to 0$, and $c_\alpha(M_{Y_n}) \to c$   for some $c \in \CC$. Then

\begin{enumerate}
\item $\displaystyle{\lim_{n\to\infty}}d(\log\ell_\alpha)(V(Y_n)) = \frac12\left(\Re\left(c^2\right) - 1\right)$
\item $\displaystyle{\lim_{n\to\infty}}d(\log\cL_\alpha)(V(Y_n)) = \frac14 c\left(\bar c^2 -1\right)$
\item $\displaystyle{\lim_{n\to\infty}}dc_\alpha(V(Y_n)) = \frac14\left(|c|^4 -2c\Re\left(c^2\right) -c^2+2c\right)$
\end{enumerate}
In particular if $Y_t$ is a flow line for $V$, $Y_n = Y_{t_n}$ for a sequence $t_n \to \infty$ and $\ell_\alpha(Y_n)\to 0$ and $c_\alpha(M_{Y_n}) \to c$ then (1)-(3) hold.
\end{theorem}

{\bf Proof:} Combining Lemmas \ref{length_der} and \ref{Fell_cont} gives (1) and (3) follows from (1) and (2).

Recalling that $V(Y_n)$ is a harmonic Beltrami differential, we will derive (2) from Lemmas \ref{three_length_der} and \ref{FL_cont} if we can show that $|J_\alpha(V(Y_n))|/\Re \cL_\alpha(M_Y) \to 0$. This will be the proven in the next section. In particular, as $\Re\cL_\alpha(M_{Y_n}) \le2 \ell_\alpha(Y_n)$ and $\ell_\alpha(Y_n) \to 0$ we have that $\Re\cL_\alpha(M_{Y_n}) \to 0$. By assumption $\Cvol(M_{Y_n})$ is bounded and $\|V(Y_n)\|_2 \to 0$. It then follows from Theorem \ref{auxillary} that $$|J_\alpha(V(Y_n))|/\Re \cL_\alpha(M_{Y_n}) \to 0$$ and (2) follows.

For the last statement we only need to show that the volume of $C(M_{Y_n})$ is bounded. As $Y_t$ is a flow line for the gradient of $-{\Rvol}$, the negative of renormalized volume, we have that the renormalized volume is bounded. However, the difference between the renormalized volume and convex core volume is bounded by a constant depending only on the the topology of the boundary and the length of the shortest compressible curve (see \cite[Theorem 1.3]{BCrenorm}). This last constant is determined by $X$ and  hence is uniform on $M_{Y_t}$. This implies that the convex core volume is bounded and the theorem follows. \eproof

\section{Bounding the norm of the auxillary term}\label{spots}
This section will be dedicated to proving the following theorem:
\begin{theorem}\label{auxillary}
There exists $\delta_0 > 0$ such that given $\eta,K>0$ there exists a $\delta >0$ such that the following holds. Assume that $Y\in \Teich(S)$ with $\Cvol(M_Y) \le K$. Let $\mu \in T_Y\Teich(S)$ be a harmonic Beltrami differential $\|\mu\|_\infty \le 3/2$ and $\alpha$ an essential closed curve on $S$ with $\|\mu\|_2 < \delta$ and $\Re\cL_\alpha(M_Y) \le \delta_0$. Then
$$|J_\alpha(\mu)| \le \eta\cdot\Re\cL_\alpha(M_Y).$$
\end{theorem}
Note that the bound of $3/2$ on $\|\mu\|_\infty$ is essentially arbitrary. We have chosen it because that is the bound that we get when $\mu$ is the harmonic Beltrami differential associated to the Schwarzian quadratic differential $\phi_Y$.

Recall that $J_\alpha(\mu)$ is the pairing of a Beltrami differential and a quadratic differential. The absolute value of the quadratic differential is a Euclidean area form and one can bound the pairing by bounding the norm of the quadratic differential and the area of the quadratic differential. Typically to show that a product of two terms is small one finds a uniform bound on one term and shows that the other term is small. To get the necessary bound here we will need to decompose $\mu$ into two parts. In one the norm will be small while in the other the support will be small while the norm will only be bounded. In particular write $\mu$ as
$$\mu = \mu^{<\epsilon}+ \mu^{\ge \epsilon}$$
where the support of $\mu^{<\epsilon}$ is the $\epsilon$-thin part of $Y$ and the support of $\mu^{\ge \epsilon}$ is the $\epsilon$-thick part of $Y$. We will then bound $J_\alpha(\mu^{<\epsilon})$ and $J_\alpha(\mu^{\ge \epsilon})$ separately.

Most of the work will be to bound $J_\alpha(\mu^{<\epsilon})$ so we begin with the easier bound on $J_\alpha(\mu^{\ge \epsilon})$.
\begin{lemma}\label{thick_bound}
Given $\epsilon, \eta>0$ there exists an $\delta >0$ such that if $\|\mu\|_2 \le \delta$ then
$$|J_\alpha(\mu^{\ge \epsilon})| \le \eta\cdot \Re \cL_\alpha(M_Y).$$
\end{lemma}

{\bf Proof:} In general bounds on $\|\mu\|_2$ don't give bounds on $\|\mu\|_\infty$. However, a mean value estimate of Teo (see \cite{Teo}) gives bounds for points in the $\epsilon$-thick part. In particular there exists a constant $C_\epsilon>0$ such that if $z \in Y^{\ge \epsilon}$ then $\|\mu(z)\| \le C_\epsilon \|\mu\|_2$. This implies
$$\left\|\left(\mu^{\ge \epsilon}\right)^{\rm aux}_{T}\right\|_\infty \leq \left\|\mu^{\ge \epsilon}\right\|_\infty \le C_\epsilon \|\mu\|_2.$$

The absolute value of the quadratic differential $dz^2/z^2$ is a Euclidean area form on $T_\alpha$. If we let $\area_\alpha$ be this area we have
$$\area_\alpha(T_\alpha) = 2\pi\cdot \Re\cL_\alpha(M_Y)$$
and therefore
$$\left|J_\alpha\left(\mu^{\ge \epsilon}\right)\right| \le \left\|\left(\mu^{\ge \epsilon}\right)^{\rm aux}_{T}\right\|_\infty \cdot  2\pi\cdot \Re\cL_\alpha(M_Y) \le 2\pi C_\epsilon \Re \cL_\alpha(M_Y) \|\mu\|_2.$$
Letting $\delta = \eta/(2\pi C_\epsilon)$ the lemma follows. \eproof

To bound $J(\mu^{<\epsilon})$ we need to bound the area of the support of $\left(\mu^{<\epsilon}\right)^{\rm aux}_T$. The proof of the following proposition will be most of the work of this section:
\begin{prop}\label{area_bound}
There exists  $\delta_0 > 0$ such that the following holds. Given $\eta,K>0$ there exists a $\epsilon>0$ such that if $\Re\cL_\alpha(M_Y) \le \delta_0$ and $\Cvol(M_Y) \le K$  
then
$$\area_\alpha\left(\supp\left(\left(\mu^{<\epsilon}\right)^{\rm aux}_T\right)\right) \le \eta\cdot \Re\cL_\alpha(M_Y).$$
\end{prop}

Assuming this for now and we can prove Theorem \ref{auxillary}.\\
\\
{\bf Proof of Theorem \ref{auxillary}:} By Proposition \ref{area_bound} we can fix an $\epsilon>0$ such that if $\Re\cL_\alpha(M_Y) \le \delta_0$ and $\Cvol(M_Y) \le K$ then
$$\area_\alpha\left(\supp\left(\left(\mu^{<\epsilon}\right)^{\rm aux}_T\right)\right) \le (\eta/3)\cdot \Re\cL_\alpha(M_Y).$$
Since $\|\mu\|_{\infty} \le 3/2$ this implies that
$$\left|J_\alpha(\mu^{<\epsilon})\right| \le (\eta/2)\cdot \Re\cL_\alpha(M_Y).$$
 By Proposition \ref{thick_bound} we can choose $\delta>0$ such that if $\|\mu\|_2< \delta$ then
$$\left|J_\alpha(\mu^{\ge \epsilon}) \right| \le (\eta/2) \cdot \Re\cL_\alpha(M_Y).$$
Note that $\delta$ depends on $\epsilon$ (and $\eta$) but $\epsilon$ only depends on $\eta$ and $K$. Therefore $\delta$ only depends on $\eta$ and $K$ and we can combine the two estimates to get the claimed bound on $|J_\alpha(\mu)|$. \eproof

The remainder of this section is dedicated to the proof of Proposition \ref{area_bound}.

\subsection{Margulis tubes}
Let $M$ be a hyperbolic n-manifold and $\epsilon >0$. Then we define the {\em thick-thin decomposition} $M = M^{\leq \epsilon}\cup M^{>\epsilon}$ with
$$M^{\leq\epsilon} = \{p\in M\ |\ \inj(p) \leq \epsilon\}\qquad M^{>\epsilon} = \{p\in M\ |\ \inj(p) >\epsilon\}.$$
By the Margulis lemma (see \cite{Thurston:book:GTTM}), there exists a constant $\epsilon_n$ such that for $\epsilon \leq \epsilon_n$ then the connected components of $M^{\leq \epsilon}$ are disjoint embedded tubular neighborhoods of cusps or simple geodesics in $M$. These tubes are called {\em Margulis tubes} and for $\alpha$ a simple closed geodesic,   we denote the tube about $\alpha$ in $M^{\leq \epsilon}$ by $T_\epsilon(\alpha)$. We further denote the radius of $T_\epsilon(\alpha)$ by $R_\epsilon(\alpha)$.

For hyperbolic 3-manifolds, if $\epsilon < \epsilon_3$ then by elementary hyperbolic geometry we have,
$$\area(\partial T_\epsilon(\alpha)) = \pi\sinh(2R_\epsilon(\alpha))\Re\mathcal L_\alpha(M).$$
As $\partial T_\epsilon(\alpha)$ has intrinsic Euclidean metric with injectivity radius $> \epsilon$ we have the inequality
$$\pi \epsilon^2 \leq \pi\sinh(2R_\epsilon(\alpha))\Re\mathcal L_\alpha(M).$$
In particular 
\begin{equation}
|\mathcal L_\alpha(M)| \geq \Re \mathcal L_\alpha(M) \geq \frac{\epsilon^2}{\sinh(2R_\epsilon(\alpha))}.
\label{tubesize}
\end{equation}

We have the following fact due to Brook-Matelski which gives uniform bounds on the change in $R_\epsilon(\alpha)$ as $\epsilon$ varies.

\begin{theorem}[{Brooks-Matelski, \cite{Brooks:Matelski:collars}}]\label{BM}
Given $\epsilon > 0$ there exists continuous functions $d^l_{\epsilon}, d^u_{\epsilon}:(0,\epsilon)\rightarrow \R_+$ such that
$d^l_{\epsilon}(\delta)\rightarrow \infty$ as $\delta\rightarrow 0$ and $d^u_{\epsilon}(\delta)\rightarrow 0$ as $\delta\rightarrow \epsilon$ and for $\alpha$ a geodesic in a hyperbolic 3-manifold with $\ell_\alpha(M) < \delta$ then
$$d^l_{\epsilon}(\delta) \leq R_\epsilon(\alpha)-R_\delta(\alpha) \leq d^u_{\epsilon}(\delta) .$$
\end{theorem}

\subsection{Orthogeodesics}
Since $Y$ is a closed surface the components of $Y^{<\epsilon}$ will be collars of closed geodesics of length $< 2\epsilon$. We let $\beta$ be one of these geodesics and $\bC_\epsilon(\beta)$ the collar. Note that it is possible that $\beta =\alpha$. The pre-image of these collars in the domain of discontinuity $\Omega_Y$ will be a collection of strips. Under the covering map from $\CC\smallsetminus\{0\}$ to
$T_\alpha$ most of these strips will map homeomorphically to strips. The only exception will be that there will be one component in the pre-image of the collar of $\alpha$ that will map to an essential annulus in $T_\alpha$. This exceptional strip will be contained in the central component $\Omega$ of $\Omega_Y$. We only need to bound the area of the strips in $T_\alpha$. We are interested in bounding the area in $T_\alpha$ of those strips that come from auxillary components. To do this we will bound the area of all strips, even those coming from the central component. Our first goal is to bound the area of a single strip. We begin with some preparation.

Associated to $\alpha$ is a solid torus cover $W_\alpha$ of $M_Y$. The surface $Y$ is a component of the conformal boundary of $M_Y$ so the $\epsilon$-collar of $\beta$ will lie in the conformal boundary of $M_Y$ and will lift to a collection of strips in the conformal boundary of $W_\alpha$. We want to bound the area of these strips in the Euclidean metric on the conformal boundary $T_\alpha = \del W_\alpha$ that comes from taking the absolute value of the quadratic differential $dz^2/z^2$. Our bounds will depend on the lengths of $\alpha$ and the constant $\epsilon$. In the hyperbolic 3-manifold the closed geodesic $\beta$ will lift to a collection of bi-infinite geodesics in $W_\alpha$ and there is a natural correspondence between the strips and these geodesics. We will see that the area of each strip decays exponentially in the distance between $\alpha$ and the corresponding lift of $\beta$.

With this informal discussion in mind we now setup the notation that we will need.  Let $\bar{M}_Y$ be the union of $M_Y$ and its conformal boundary. Then the cover $W_\alpha$ extends to a cover $\bar W_\alpha$ of $\bar{M}_Y$. Let
$$\pi_\alpha\colon \bar W_\alpha \to \bar M_Y$$
be the covering map. We also have the nearest point projection
$$r\colon M_Y \to C(M_Y)$$
to the convex core. This map extends continuously to $\bar M_Y$ and lifts to a map
$$r_\alpha\colon \bar W_\alpha\to C(W_\alpha)$$
where $C(W_\alpha) = \pi^{-1}_\alpha(C(M_Y))$. In particular we have $r\circ\pi_\alpha = \pi_\alpha \circ r_\alpha$. 

\begin{figure}[htbp] 
   \centering
   \includegraphics[width=3in]{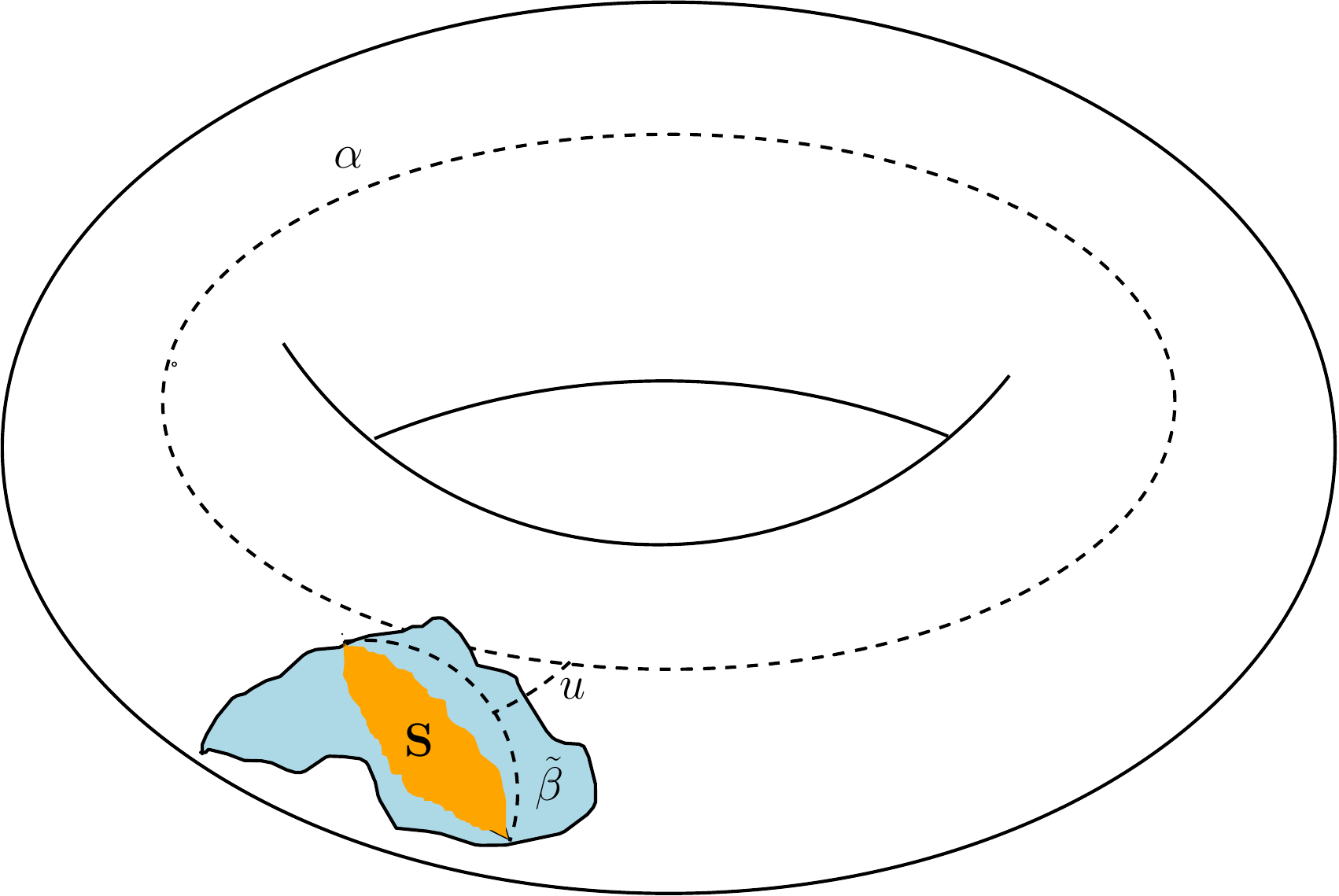} 
   \caption{Component $\bS$ of $S^\alpha_\epsilon(\beta)$}
   \label{Sab}
\end{figure}

Now let $S^\alpha_\epsilon(\beta)$ be the collection of strips in $\del\bar W_\alpha$ that map to $\bC_\epsilon(\beta)$. If $\beta\neq \alpha$ then this is just the union of components of $\pi_\alpha^{-1}(\bC_\epsilon(\beta))$. If $\beta=\alpha$ then there is one annular component in $\pi_\alpha^{-1}(\bC_\epsilon(\alpha))$ that is not included in $S^\alpha_\epsilon(\alpha)$. By \cite{EMM1} the retraction $r$ is a $2$-Lipschitz map from the hyperbolic metric on $Y$ to induced path metric on $\del C(M_Y)$. Therefore $r(\bC_\epsilon(\beta)) \subset T_{2\epsilon}(\beta)$ and if $\bS$ is a component of $S^\alpha_\epsilon(\beta)$ there is a unique component $\bT$ of $\pi_\alpha^{-1}(T_{2\epsilon}(\beta))$ with $r_\alpha(\bS) \subset \bT$ (see figure \ref{Sab}).

There is a unique shortest geodesic $u$ in $W_\alpha$ from $T_{\epsilon_3}(\alpha)$ to $\bT$. This geodesic $u$ will be orthogonal to the boundary of the two tubes so we call $u$ an {\em orthogeodesic}. We will see that the area of $\bS$ decays exponentially in the length of $u$. We note that while it may seem more natural to take the orthogeodesic from the gedoesic $\alpha$ to the geodesic in the core of $\bT$, for our purposes our choice of orthogeodesic is more convenient. We also emphasize that one endpoint of our orthogeodesic will always lie on the boundary of the $\epsilon_3$-Margulis tube $T_{\epsilon_3}(\alpha)$ while the other end will lie on a lift of the $2\epsilon$-Margulis tube $T_{2\epsilon}(\beta)$ so while our notation does not emphasize this, the orthogeodesic $u$ depends on $\epsilon$ and as $\epsilon$ decreases the orthogeodesic will become longer.

If $H$ is an embedded half-space in $W_\alpha$ then its closure in the conformal boundary is a round disk $D$. We begin by bounding the area of $D$ in terms of the distance between the core geodesic of $W_\alpha$ and $H$.

\begin{lemma}\label{area-disk}
Let $H$ be an isometrically immersed half-space in $W_\alpha$ that is disjoint from the core geodesic $\alpha$ with boundary disk $D$. Then
$$\area_\alpha(D) \leq  \frac{\pi}{\sinh^2(d(\alpha,H))}.$$
\end{lemma}
{\bf Proof:} We can lift the picture to $\Hs$ and assume that $\alpha$ has endpoints  $\pm 1$ and $D$ is the disk of radius $r$ centered at $0$ with $d(\alpha, H) = e^{-r}$. Note that while the disk $D$ is embedded in $\CC$ it may be be immersed in $\del W_\alpha$. However, this will only decrease its area. The map $f(z)= (z-1)/(z+1)$ sends the geodesic $\alpha$ to the geodesic with endpoints $0$ and $\infty$. Thus pulling back the metric $|dz|^2/|z|^2$ we have
$$d\area_\alpha = \frac{|f'(z)|^2}{|f(z)|^2} |dz|^2 = \frac{4}{|z-1|^2|z+1|^2} |dz|^2 =\frac{4}{|z^2-1|^2} |dz|^2.$$
Therefore on $D$ we have 
$$d\area_\alpha \le \frac4{|r^2 - 1|^2} |dz|^2$$
and it follows that
$$\area_\alpha(D) \leq \frac{4 \pi r^2}{(1-r^2)^2}  =  \frac{\pi}{\sinh^2(d(\alpha,H))}.$$
\eproof

In the next lemma we show that each component of $S^\alpha_\epsilon(\beta)$ lies in a round disk in $T_\alpha = \del W_\alpha$ which bounds a half-space in $W_\alpha$ whose distance from the core Margulis tube $T_{\epsilon_3}(\alpha)$ is given in terms of the length of the associated orthogeodesic.
\begin{lemma}\label{shadow}
There exists a $\delta_0 > 0$ such that if $\Re \cL_\alpha(M_Y) <\delta_0$ and $\epsilon < \epsilon_3/2$ then the following holds. Let $\ell_\beta(Y) < 2\epsilon$ and $\bS$ be a component of $S^\alpha_\epsilon(\beta)$.  Let $\bT$ be the component of $\pi_\alpha^{-1}(T_{2\epsilon}(\beta))$ with $r_\alpha(\bS) \subset \bT$ and let $u$ be the orthogeodesic between $T_{\epsilon_3}(\alpha)$ and $\bT$.
Then there is a half-space $H$ in $W_\alpha$ whose boundary disk $D$ contains $\bS$ with 
$$d(\alpha, H) = R_{\epsilon_3}(\alpha)+ \ell(u) - \log\left(1+\sqrt2\right) \ge \log\sqrt 2.$$
\end{lemma}
{\bf Proof:} 
By the inequality \ref{tubesize} for the radius of the Margulis tube we can choose $\delta_0 > 0$ such that if $\Re \cL_\alpha(M_Y) <\delta_0$ then
$$R_{\epsilon_3}(\alpha) \geq \log(1+\sqrt{2}) + \sqrt{2}.$$
(Choosing $\delta_0 = \epsilon_3/50$ will suffice). We also note that the inequality in the statement of the Lemma follows from the lower bound on $R_{\epsilon_3}(\alpha)$.

As we assume $\ell_\beta(Y) < 2\epsilon$, it follows that $S^\alpha_\epsilon(\beta)$ is non-empty so we can consider a component $\bS$. We again lift the picture to $\Hs$. Let $\tilde\bT$ be a lift of $\bT$ to $\Hs$ and let $\tilde u$ be a lift of $u$ that is orthogonal to $\del\tilde \bT$. Further we extend $u$ to the perpendicular $v$ from $\partial\bT$ to $\alpha$ with lift $\tilde v$ containing $\tilde u$. We can assume that $\tilde u$ is the vertical geodesic segment $\{0\} \times [1, e^{\ell(u)}]$ in the upper-half-space model $\Hs = \CC \times \R_+$ and that $\tilde\bT$ is orthogonal to $\tilde u$ at $(0, 1)$. Therefore  $\tilde v$ is the geodesic $\{0\}\times [1,e^{\ell(v)}]$ where $\ell(v) = R_{\epsilon_3}(\alpha) +\ell(u)$.  Let $H'$ be the half-space in $\Hs$ that contains $\tilde\bT$ and whose boundary plane is orthogonal to $\tilde u$ at $(0,1) \in \Hs$.

The retraction $r_\alpha$ lifts to a retraction
$$\tilde r\colon {\bar{\mathbb H}}^3 \to C(\Lambda)$$
where $\Lambda$ is the limit set of the Kleinian group uniformizing $M_Y$. Let $\tilde\bS$ be the lift of $\bS$ with $\tilde r\left(\tilde\bS\right) \subset \tilde\bT$. Let $\tilde D \subset \CC$ be the Euclidean disk of radius $1+ \sqrt 2$ centered at $0 \in \CC$ and $\tilde H$ the half-space in $\Hs$ with boundary $\tilde D$. We'll show that $\tilde\bS\subset \tilde D$ so that if $D$ is the image of $\tilde D$ in $\del W_\alpha$ then $\bS\subset D$. The immersed half-space $H \subset W_\alpha$ bounded by $D$ is the image of $\tilde H$ so we have that $d(\alpha, H)$ is equal to $\ell(v)$ minus the distance from  from the boundary plane of $\tilde H$ to $H'$. As this latter distance is $\log\left(1+\sqrt2\right)$  this will give
$$d(\alpha, H) = R_{\epsilon_3}(\alpha) +\ell(u) - \log\left(1+\sqrt2\right).$$
If we can show the inclusion $\tilde\bS \subset \tilde D$ we are done.

\begin{figure}[htbp] 
   \centering
   \includegraphics[width=3in]{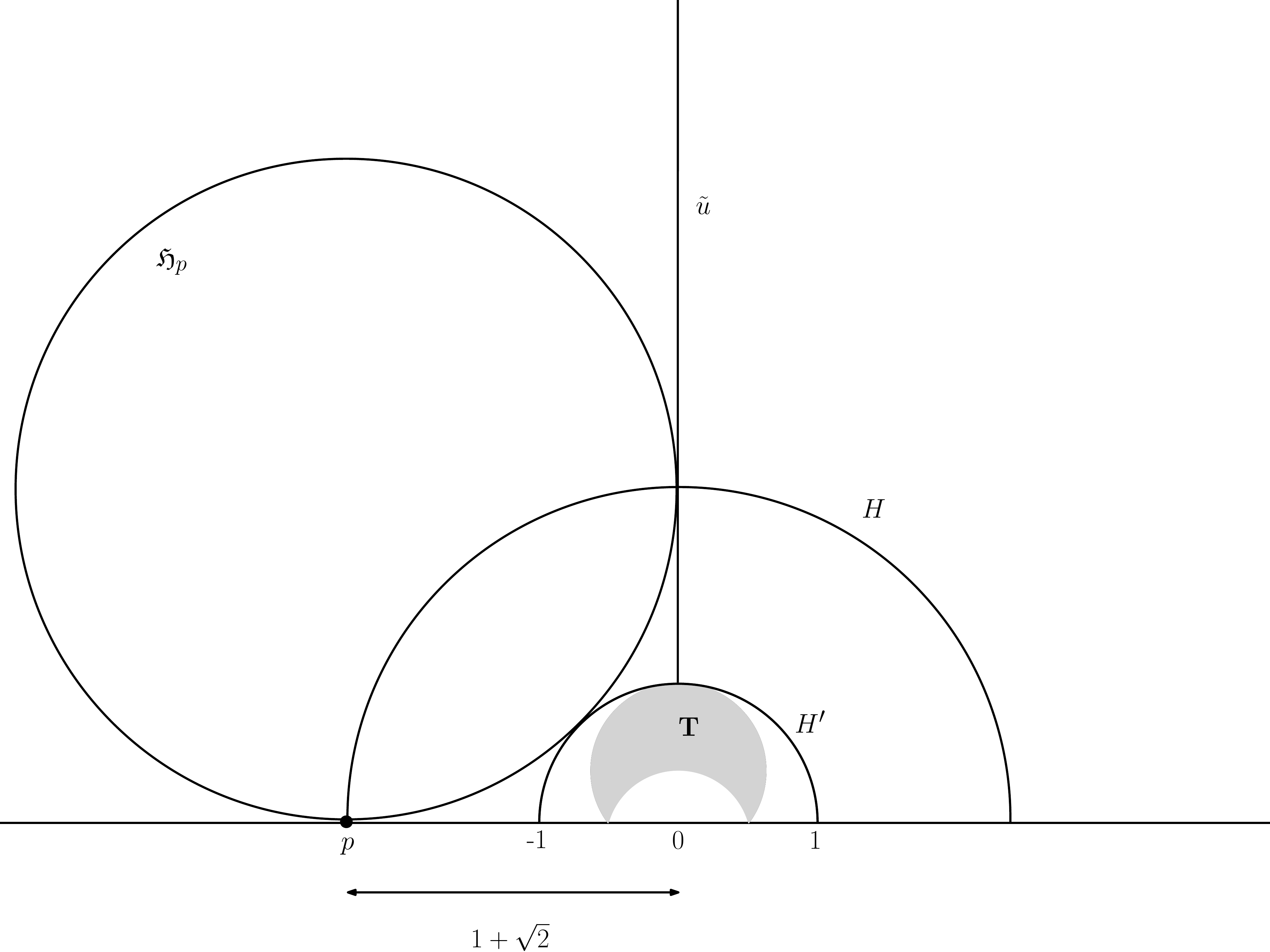} 
   \caption{Tubes, half-spaces and horoballs}
   \label{fig:tubes}
\end{figure}

We now prove this inclusion. For $z\in \chat\smallsetminus \Lambda$ there is a unique horosphere $\mathfrak H_z$ that intersects $C(\Lambda)$ in a single point. This point of intersection is $\tilde r(z)$. If $\tilde r(z) \in \tilde\bT$
then $\mathfrak H_z$ will intersect $H'$. The perpendicular $\tilde v $   from $\tilde\alpha$ to $\tilde\bT$ is contained in $C(\Lambda)$ and therefore must not  
intersect the interior of $\mathfrak H_z$. As $\tilde v$ is the vertical geodesic $\{0\}\times [1,e^{\ell(v)}]$ a simple calculation shows that when $\ell(v) \geq \log(1+\sqrt{2})$ the interior of $\mathfrak H_z$ will intersect $\tilde v$ if $|z|> 1+\sqrt2$. (See Figure \ref{fig:tubes}.) By our choice of $\delta_0$, we have $\ell(v) = R_{\epsilon_3}(\alpha)+\ell(u) \geq \log(1+\sqrt{2})$. Therefore if $z\in \tilde\bS$ we have that $z\in \tilde D$ and $\tilde\bS \subset\tilde D$, completing the proof.
\eproof
\eproof

Let $\cO_\epsilon^\alpha(\beta)$ be the set orthogeodesics from $T_{\epsilon_3}(\alpha)$ to the components of $\pi_\alpha^{-1}(T_{\epsilon}(\beta))$. With $\delta_0$ as in the above Lemma \ref{shadow}, we have.

\begin{corollary}\label{thin-area}
 Given $\delta, \eta>0$   there exists $\epsilon>0$ such that if $\Re\cL_\alpha(M_Y) < \delta_0$ and $\Re\cL_\beta(M_Y) < 2\delta$ then
 $$\area_\alpha(S^\alpha_\epsilon(\beta)) \leq \eta\cdot\Re \cL_\alpha(M_Y) \sum_{u\in \cO^\alpha_\delta(\beta)} e^{-2\ell(u)}.$$
\end{corollary}

{\bf Proof:} 
We first show that for $\epsilon < \epsilon_3/2$  and $\ell_\beta(Y) \leq 2\epsilon$ then 
$$\area_\alpha(S^\alpha_\epsilon(\beta)) \leq A\cdot\Re \cL_\alpha(M_Y)\sum_{u\in \cO^\alpha_{2\epsilon}(\beta)} e^{-2\ell(u)}$$
where $A$ is a universal constant.

 As $\ell_\beta(Y) \leq 2\epsilon$ then $S^\alpha_\epsilon(\beta)$ is non-empty. Given a component $\bS$ of $S^\alpha_\epsilon(\beta)$ let $\bT$ be the component of $\pi_\alpha^{-1}(T_{2\epsilon}(\beta))$ with $r_\alpha(\bS) \subset \bT$ and let $u$ be the orthogeodesic between $T_{\epsilon_3}(\alpha)$ and $\bT$. Then as $\Re\cL_\alpha(M_Y) < \delta_0$, by Lemma \ref{shadow}, there is an immersed half-space $H$ in $W_\alpha$ whose boundary disk $D$ contains $\bS$ with
\begin{eqnarray*}
d(\alpha, H) 
& = & R_{\epsilon_3}(\alpha) + \ell(u) - \log\left(1+\sqrt2\right).
\end{eqnarray*}
If $t\ge \log\sqrt2$ then $\sinh(t) \ge e^t/4$. Since the disk $D$ contains $\bS$ and by Lemma \ref{shadow} we have $d(\alpha, H) \ge \sqrt{2}$ we can combine this estimate on $\sinh$ with Lemma \ref{area-disk} to get
$$\area_\alpha(\bS) \le 16\pi e^{-2d(\alpha, H)} \le 16\pi e^{-2\left(R_{\epsilon_3}(\alpha) + \ell(u) -\log\left(1+\sqrt2\right)\right)}.$$ 

By \eqref{tubesize} we have
$$\Re\cL_\alpha(M_Y) \ge \frac{\epsilon_3^2}{\sinh(2R_{\epsilon_3}(\alpha))} \ge 2\epsilon_3^2 e^{-2R_{\epsilon_3}(\alpha)}.$$

Together the two estimates give
$$\area_\alpha(\bS) \le \left(8\pi/\epsilon_3^2\right)\cdot\left(3+2\sqrt2\right)\cdot \Re \cL_\alpha(M_Y) e^{-2\ell(u)}.$$
Summing over the set of all orthogeodesics in $\cO^\alpha_{2\epsilon}(\beta)$ gives  our first estimate.

Now let $2\epsilon < \delta$. We let $\bT'$ be the component of $\pi^{-1}_\alpha\left(T_{\delta}(\beta)\right)$ that contains $\bT$ and let $u'$ be the orthogeodesic between $T_{\epsilon_3}(\alpha)$ and $\bT'$. Note that  the map $u\mapsto u'$ defines an bijection between $\cO^\alpha_{2\epsilon}(\beta)$ and $\cO^\alpha_{\delta}(\beta)$. Then $u'$ is a subsegment of $u$ and
$$\ell(u) -\ell(u') = R_{\delta}(\beta) - R_{2\epsilon}(\beta) \ge d^l_{\delta}(2\epsilon)$$
where the function on the right comes from Theorem \ref{BM}. 
Thus for $\ell_\beta(Y) \leq 2\epsilon$
$$\area_\alpha(S^\alpha_\epsilon(\beta)) \leq Ae^{-2d^l_{\delta}(2\epsilon)}\cdot\Re \cL_\alpha(M_Y) \sum_{u\in \cO^\alpha_\delta(\beta)} e^{-2\ell(u)}.$$
For $\ell_\beta(Y) > 2\epsilon$ the inequality also holds trivially as $S^\alpha_\epsilon(\beta) = \emptyset$ and as $\Re\cL_\beta(M_Y) <2\delta$, the righthandside is positive.
Thus as $ d^l_{\delta}(2\epsilon)\to \infty$ as $\epsilon \to 0$ the theorem follows.
\eproof

\subsection{Uniform bounds on the Poincar\'e series}
Except for the last lemma, in this subsection $M= \Hs/\Gamma$ can be any complete hyperbolic 3-manifold uniformized by a Kleinian group $\Gamma$. 

For $x\in \Hs$ we define the {\em Poincar\'e series}
$$P_\alpha(x) = \sum_{\gamma\in\Gamma} e^{-2d(x, \gamma x)}$$
where $d$ is the hyperbolic distance. Poincar\'e series play an important role in dynamics but here we will only be interested in obtaining uniform bounds on $P_\alpha(x)$ which we will then use to bound the similar exponential sum of lengths of orthogeodesics that appears in Corollary \ref{thin-area}.

We have the following elementary bound.

\begin{lemma}\label{proj_est}
\phantom{ }
\begin{itemize}
\item Let $\calC$ be a closed convex set in $\Hs$ and $\bar\calC$ is closure in $\barHs =\Hs \cup \chat$ and let
$$r\colon \barHs\to \bar\calC$$
be the nearest point projection. If $x \in \del\calC$ and $B_x(\epsilon)$ is the $\epsilon$ ball centered at $x$ then $r^{-1}(B_x(\epsilon))$ contains a round disk $D \subset \chat$ that bounds a half-space $H\subset \Hs$ with 
$$e^{d(x,H)} \le 2\coth(\epsilon).$$

\item Given $x \in \Hs$ and a half-space $H\subset \Hs$ with boundary disk $D \subset \chat$ we have
$$\area_x(D) \ge \pi e^{-2d(x,H)}$$
where $\area_x$ is the visual measure on $\chat$ determined by $x$.

\end{itemize}
\end{lemma}

{\bf Proof:} Let $P$ be a support plane for $\calC$ at $x$. Then $P$ bounds a half-space whose interior is disjoint from the interior of $\calC$ and we let $D$ be the round disk in the conformal boundary of this half-space that projects orthogonally to $P$ with image in $P \cap B_x(\epsilon)$. Let $y$ be the point on $H$ closest to $x$, $z$ a point in $\del D$, and $w$ the orthogonal projection of $z$ to $P$. Then $wxyz$ forms a planar quadrilateral with one ideal vertex at $z$ and right angles at all other vertices. Standard formulas give $\sinh d(x,y) \sinh d(x,w) = 1$ which can be rewritten in a less symmetric form as $\cosh d(x,y) = \coth d(x,w)$. Noting that $d(x,y) = d(x,H)$, $d(x,w) = \epsilon$ and $e^t \le 2\cosh t$ we have $e^{d(x,H)} \le 2 \coth(\epsilon)$ so we will be done if we can show that $r(D) \subset B_x(\epsilon)$.

For $z \in D$ let $\gamma$ be the geodesic in $\barHs$ from $z$ to $r(z)$ and let $p = P \cap \gamma$. Note that $r(z) = r(p)$ and since $r$ is a contraction (and $r(x) = x$) we have that $d(x, r(z)) \le d(x, p)$. Let $B'$ be the ball of radius $d(x,p)$ centered at $x$ and let $H_p$ be the closed half-space whose boundary is tangent to $B'$ at $p$ and whose interior is disjoint from $x$. If $r(z) \neq p$  then $r(z)$ is in the interior of $B'$ and is hence disjoint from $H_p$. As $\gamma$ intersects the boundary plane of $H_p$ we have that $z$, the other endpoint of $\gamma$, will be contained in $H_p$. As the boundary plane of $H_p$ is orthogonal to $P$,  the orthogonal projection of every point in $H_p$ to $P$ will have distance from $x$ that is at least $d(x,p)$ so we have $d(x,p) < \epsilon$. Combining inequalities gives $d(x, r(z)) \le \epsilon$ when $r(z) \neq p$. If $r(z) = p$ then the interior of $\gamma$ will be disjoint from $B'$ as for every $q \in \gamma$ we have $r(q) = r(z) = p$ and $d(q,x) \ge d(p,x)$ by the contraction of $r$. However, if $\gamma$ is not contained in $H_p$ it must intersect the interior of $B'$. Therefore $\gamma$, and in particular $z$, lie in $H_p$. As above this implies that 
$r(z) \in B_x(\epsilon)$ completing the proof that $r(D) \subset B_x(\epsilon)$.

\begin{figure}[htbp] 
   \centering
\includegraphics[width=2.5in]{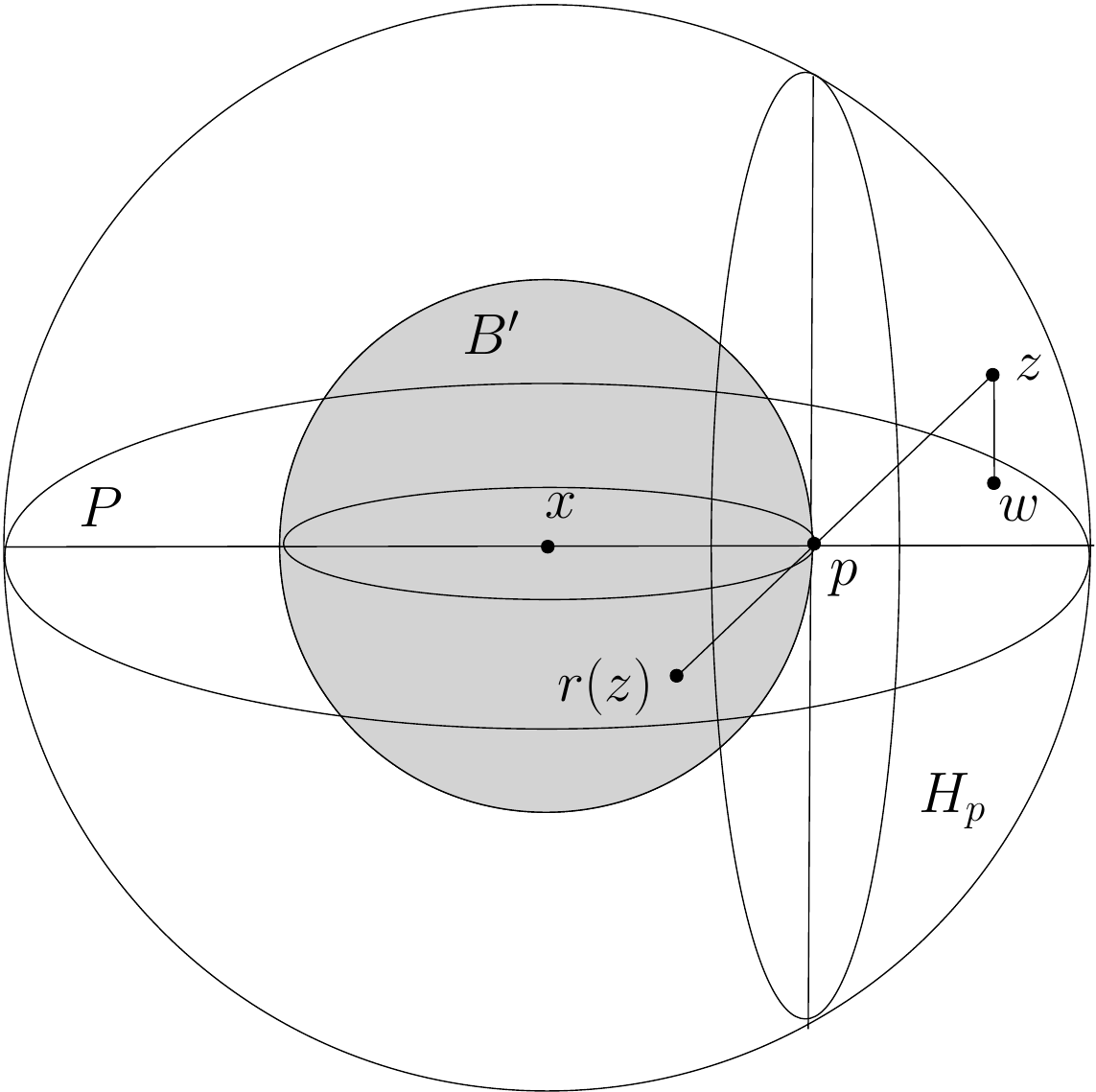} 
   \caption{View in Klein model}
   \label{klein}
\end{figure}

For the second bullet, as above let $y$ be the point on $H$ closest to $x$ and let $z$ be a point in $\del D$. These three points from a right triangle with one idea angle at $z$ and angle $\theta$ at $x$. Therefore $D$ is a spherical disk of radius $\theta$ in the visual metric based at $x$. For such right triangles we have $\sech d(x,y) = \sin\theta$. The area of spherical triangle is $2\pi(1-\cos\theta)$ so $\sin^2\theta = 1-\cos^2\theta \le \area_x(D)/\pi$. Since $e^{-2d(x,y)} \le \sech d(x,y)$ and $d(x,y) = d(x,H)$ we can combine expressions to get
$$\pi e^{-2d(x,H)} \le \area_x(D).$$
\eproof

Using these estimates we can obtain our uniform bound on the Poincar\'e series.
\begin{lemma}\label{unif-poinc}
If $M =\Hs/\Gamma$ is a hyperbolic 3-manifold and $x$ a point on the boundary of the convex core $C(M)$ with injectivity radius $\geq \epsilon$, then for $\tilde x$  a lift of $x$ to $\Hs$ we have
$$P_\alpha(\tilde x) \leq 8\coth^2(\epsilon).$$
\end{lemma}

{\bf Proof:} Let $\Lambda$ be the limit set of $\Gamma$ and  $r:\barHs \rightarrow C(\Lambda)$ the retract map to the convex hull of the limit set. Let $B_\epsilon(\tilde x)$ be the ball of radius $\epsilon$ about $\tilde x $. As the injecitivity radius of $x$ is $\ge \epsilon$ we have $B_\epsilon \cap \gamma(B_\epsilon) = \emptyset$ for all $\gamma \in \Gamma\smallsetminus\{\id\}$. We define $U = r^{-1}(B_\epsilon)\cap \chat$. Since $r$ is $\Gamma$-invariant we have $U \cap \gamma(U) = r^{-1}(B_\epsilon\cap \gamma(B_\epsilon))= \emptyset$ if $\gamma \in \Gamma\smallsetminus \{\id\}$. By Lemma \ref{proj_est} we have that $U$ contains a disk $D$ that bounds a half-space $H$ with $e^{d(\tilde x, H)} \le 2 \coth(\epsilon)$. By the triangle inequality for all $\gamma \in \Gamma$ we have that $d(\tilde x,\gamma(H)) \le d(\tilde x, \gamma(\tilde x)) +d(\gamma(\tilde x), \gamma(H))$. As $d(\gamma(\tilde x), \gamma(H)) = d(\tilde x, H)$ this gives $d(\tilde x, \gamma(\tilde x)) \ge d(\tilde x, \gamma(H)) - d(\tilde x, H)$ which in turn implies
\begin{eqnarray*}
\sum_{\gamma\in \Gamma} e^{-2d(\tilde x, \gamma(\tilde x))} & \le &e^{2d(\tilde x, H)} \sum_{\gamma} e^{-2d(\tilde x, \gamma(H))}\\
& \le & 2\pi^{-1}\coth^2(\epsilon) \sum_{\gamma\in \Gamma} \area_{\tilde x}(\gamma(D))\\
&\le & 8\coth^2(\epsilon)
\end{eqnarray*}
where the last sum is bounded by $\area_{\tilde x}(\chat) = 4\pi$ since the disks $\gamma(D)$ are all disjoint.
\eproof

We now can obtain a uniform bound on the exponential sum of orthogeodesic lengths.
\begin{lemma}\label{ortho2poincare}
Fix $\delta<\epsilon_3$ and assume that $\alpha$ and $\beta$ are closed geodesics in $M_Y \in CC(N,S;X)$ with length $\le 2\delta$. Then the sum
$$\sum_{u\in \cO^\alpha_\delta(\beta)} e^{-2\ell(u)}$$
is bounded by a constant that only depends on $\delta$ and the diameter of $C(M_Y)^{>\delta}$.
\end{lemma}

{\bf Proof:} Choose $x \in \del C(M_Y) \cap C(M_Y)^{\ge \epsilon_3}$ and let $\tilde x \in \Hs$ be in the pre-image of $x$ under the covering map $\Hs \to M_Y$. We will define a map $\psi$ from $\cO^\alpha_\delta(\beta)$ to $\Gamma$ such that $\ell(u) \ge d(\tilde x, \psi(u)(\tilde x)) - 2D$  where $D$ is the diameter of $C(M_Y)^{>\delta}$ and $\psi$ is at most $N$-to-one. Once we have defined such a $\psi$ we have that
$$\sum_{u\in \cO^\alpha_\delta(\beta)} e^{-2\ell(u)} \le \sum_{u\in \cO^\alpha_\delta(\beta)} e^{-2d(\tilde x, \psi(u)(\tilde x))+2D} \le Ne^{2D} \sum_{\gamma\in\Gamma} e^{-2d(\tilde x, \gamma(\tilde x))}$$

We claim that there exists an integer $N_0$ such that at most $N_0$ components of the pre-image of $T_{\epsilon_3}(\alpha)$ and $T_{\epsilon_3}(\beta)$ in $\Hs$  intersect any ball of radius $D$.  To see this we note that as $\beta$ has length less than $2\delta$, $T_{\delta}(\beta)$ is non-empty. Thus each point in $T_{\epsilon_3}(\beta)$ is contained in a ball of radius $r = d^l_{\epsilon_3}(\delta)$ which lies inside $T_{\epsilon_3}(\beta)$. Therefore there are $N_0$ disjoint balls of radius $r$ in the ball of radius $D+2r$. Thus $N_0\vol(B_r) \leq\vol(B_{D+2r})$ bounding $N_0$. Let $N = N_0^2$.

For each $u \in \cO^\alpha_\delta(\beta)$ we will carefully choose a lift $\tilde u$ to $\Hs$. More precisely we will carefully choose where the initial endpoint $\tilde u_-$ lies and then the terminal endpoint $\tilde u_+$ will determine $\psi(u)$. Note that the endpoints $\tilde u_-$ and $\tilde u_+$ will lie in the pre-image of $C(M)^{>\delta}$ so we can choose $\tilde u_-$ such that $d(\tilde x, \tilde u_-) \le D$. Then there will be a $\gamma\in \Gamma$ with $d(\gamma(\tilde x), \tilde u_+)\le D$ and we define $\psi(u) = \gamma$ and it follows that
$$\ell(u) \ge d(\tilde x, \psi(u)(\tilde x)) - 2D.$$

The orthogeodesic $u$ is determined by which components  in the pre-image of $T_{\epsilon_3}(\alpha)$ and $T_{\epsilon_3}(\beta)$  that $\tilde u_-$ and $\tilde u_+$ lie in. As there are at most $N_0$ in a $D$-neighborhood of $\tilde x$ and $N_0$ in a $D$-neighborhood of $\gamma(\tilde x)$ there are at most $N = N_0^2$ orthogeodesics $u$ with $\psi(u) = N$ so $\psi$ is at most $N$-to-one. \eproof

For the previous lemma to be useful need to able to control the diameter of $C(M_Y)^{>\epsilon}$. The following argument is well known.
\begin{lemma}
Let $\epsilon > 0, V> 0$ and $M$ a geometrically finite hyperbolic manifold  and convex core volume less than $V$. Then the diameter of every connected component of  $C(M)^{>\epsilon}$, the $\epsilon$-thick part of the convex core, is bounded by a constant that is a function of $\epsilon, V, \chi(\partial M)$ and length of the shortest compressible curve in $\del C(M)$.\label{diam}
\end{lemma}
{\bf Proof:}
 Given $x,y$ in connected component of $C(M)^{> \epsilon}$ let $\gamma:[0,1]\rightarrow  C(M)^{> \epsilon}$ be a path joining $x,y$. Then $f(t) = d(x, \gamma(t))$ is continuous with $f(0) = 0$ and $f(1) = d(x,y)$ and we can take $t_k \in [0,1]$ with $f(t_k) = 2k\epsilon$ for $0\le k \leq n$ where $2n\epsilon \leq d(x,y) < 2(n+1)\epsilon$. Then the balls $B(\gamma(t_k), \epsilon)$ are disjoint, embedded and contained in  $N_\epsilon(C(M))$. Thus the number of disjoint $\epsilon$-balls is bounded in terms of the volume of $N_\epsilon(C(M))$ which in turn bounds $d(x,y)$. 

To bound the volume of $N_\epsilon(C(M))$, let $V_t$ be the volume of $N_t(C(M))$, then $\dot V_t = A_t$ the area of $\partial N_t(C(M))$. If $\beta_M$ is the bending lamination on $\partial C(M)$ then an easy calculation gives
$$A_t = 2\pi|\chi(\partial M)|\cosh^2(t) +L(\beta_M)\sinh(t)\cosh(t)$$
where $L(\beta_M)$ is the length of $\beta_M$. Also by \cite{BBB}, there are universal constants $A, B$ such that $L(\beta_M) \leq (A + B/\delta)|\chi(\partial M)| $. Integrating, it follows that the volume of $N_\epsilon(C(M))$ is bounded by a function of  $V$, $\chi(\del M)$ and $\delta$.
\eproof

To apply the previous lemma to bound the diameter of $C(M)^{>\epsilon}$ we need to know that this set is connected. In general this will be false. The next lemma show that in our setting it holds for sufficiently small $\epsilon$ for hyperbolic manifolds in a relatively acylindrical deformation space.
\begin{lemma}\label{unif_conn}
Given $CC(N,S;X)$ with $(N,S)$  relatively acylindrical then there exists a  $\delta > 0$  such that if $M \in CC(N,S;X)$ then $C(M)^{>\delta}$ is connected.
\end{lemma}

{\bf Proof:} 
Let $X'$ be the union of components of $\del C(M)$ that face $X$.  The retraction from $X$ to $\del C(M)$ is Lipschitz with Lipschitz constant only depending on $\inj(X)$ (see  \cite{BC_retract_lip}). Therefore the  diameter of $X$ is bounded by $D$, a constant depending only on $\chi(X)$ and $\inj(X)$.

Let $Y'$ be the union of components of $\del C(M)$ facing $Y$ and recall that the path metric on $Y'$ is hyperbolic. Since $(N;S)$ is relatively  acylindrical there is at most one homotopy class of curve in $S$ that is homotopic to $\alpha$. If there is such a curve, let $C$ be the $\epsilon_3$-Margulis collar about its geodesic representative $\alpha'$ in $Y'$. Note that if there is no such curve or the curve has length $\ge 2\epsilon_3$ then $C$ is empty. By a theorem of Bers $Y'$ has a bounded length pants decomposition with constants only depending on the topology of $Y'$ and each point in the $\epsilon_3$-thick part of $Y'$ is a uniformly bounded distance from at least two of these curves. On the other hand if $p$ is in the $\epsilon_3$-thin part but not in $C$ there is an essential curve of length $\le 2\epsilon_3$ that is not homotopic to a multiple of $\alpha$. Therefore there exists an $L > 0$ such that for $p \in Y'-C$ there is a closed curve of length $\leq L$ containing $p$ which is not homotopic to a multiple of $\alpha$.

We now choose $\delta > 0 $ such that $d^l_{\epsilon_3}(\delta) \geq \max\{D,L\}$. Let $p$ be a point in $\del C(M) \cap T_\delta(\alpha)$. If $p \in X'$ then the diameter bound on $X'$ implies that $X' \subset T_{\epsilon_3}(\alpha)$. This is a contradiction since $X'$ cannot be contained in a Margulis tube. If $p \in Y' \smallsetminus C$ then there is an essential closed curve $\beta$ through $p$ of length $\le L$ that is not homotopic to $\alpha$. Then the length bound implies that $\beta \subset T_{\epsilon_3}(\alpha)$ which is again a contradiction. It follows that $T_{\delta}(\alpha) \cap \partial C(M) \subseteq C$.

To finish the proof we show that $\del T_{\delta}(\alpha)\cap C(M)$ is connected. We can assume that $\del T_\delta(\alpha)$ and $\del C(M)$ are transverse. (If not we can slightly decrease $\delta$.) On $\del T_\delta(\alpha)$ the intersection $\del T_\delta(\alpha) \cap \del C(M)$ will be a collection of simple closed curves that are either homotopic to $\alpha$ or are contractible and bound disks. These curves will bound the region of $\del T_\delta(\alpha)$ that is contained in $C(M)$. We'll show that it is either a pair of parallel curves homotopic to $\alpha$ or a collection of contractible curves bounding disks whose interiors are disjoint from $C(M)$.

First we assume that  the intersection of $\alpha'$ with the bending lamination is zero. Then the collar $C$ is totally geodesic outside of $\alpha'$ where it is bent at some angle (possibly zero). This implies that for points in $C$ the injectivity radius on $\del C(M)$ agrees with the injectivity radius in the ambient hyperbolic 3-manifold $M$ and therefore the intersection of $\del C(M)$ with $T_{\delta}(\alpha)$ is the $\delta$-Margulis collar about $\alpha$. In particular $\del T_\delta(\alpha) \cap \del C(M)$ is a pair of parallel curves and  $\del T_\delta(\alpha) \cap C(M)$ is an annulus bounded by these curves  and therefore is connected.
If $\alpha'$ intersects the bending lamination then we can foliate $C$ with geodesic segments in  $M$ joining the boundary components of $C$. Note that $C$ is contained in $T_{\epsilon_3}(\alpha)$ and as any geodesic in $T_{\epsilon_3}(\alpha)$ will intersect $T_\delta(\alpha)$ in a connected set we have that each geodesic arc of the foliation has connected intersection with $T_\delta(\alpha)$. Thus the intersection is either empty, a point or a closed interval. By continuity of the foliation this implies that $C\cap T_{\delta}(\alpha)$ is either a union of disks or an annulus. In the annulus case, then  $\partial T_\delta(\alpha)\cap C(M)$ is also an annulus and therefore is connected. For the disk case, $Y'$ intersects $T_\delta(\alpha)$ in a union of disjoint disks and $\partial T_\delta(\alpha)\cap C(M)$ is either a union of disks or the complement of a union of disks. As $T_\delta(\alpha)\cap C(M)$ is the intersection of two convex sets, it is connected. Therefore $\partial T_\delta(\alpha)\cap C(M)$ must be the complement of a union of disks and therefore is connected.
\eproof

\subsection{Proof of Proposition \ref{area_bound}}
We now are ready to prove Proposition \ref{area_bound} which will complete the proof that the auxillary components do not contribute to the limiting model flow. We restate it.\newline

\noindent
{\bf Proposition \ref{area_bound}}
{\em There exists  $\delta > 0$ such that the following holds. Given $\eta,K>0$ there exists a $\epsilon>0$ such that if $\Re\cL_\alpha(M_Y) \le 2\delta$ and $\Cvol(M_Y) \le K$ 
then
$$\area_\alpha\left(\supp\left(\left(\mu^{<\epsilon}\right)^{\rm aux}_{T}\right)\right) \le \eta\cdot \Re\cL_\alpha(M_Y).$$
}

{\bf Proof:} The support of $\left(\mu^{<\epsilon}\right)^{\rm aux}_{T}$ is contained in the union of $S_\epsilon^\alpha(\beta)$ where $\beta$ ranges over the closed geodesics of length $<2\epsilon$ in $Y$. There are at most $3g-3$ such curves, where $g$ is the genus of $Y$, so it will be enough to bound the area of the individual $S_\epsilon^\alpha(\beta)$.

Fix $\delta>0$ as in Lemma \ref{unif_conn} and
choose $\tilde x \in \Hs$ that maps to $x \in M^{\ge \epsilon_3} \cap \del C(M)$ under the covering map $\Hs \to M_Y$. Then by Lemmas \ref{unif-poinc} and \ref{ortho2poincare} we  have
for $\alpha,\beta$ of length less than $2\delta$ $$\sum_{u\in \cO^\alpha_\delta(\beta)} e^{-2\ell(u)} \le Ne^{2D}\sum_{\gamma\in \Gamma} e^{-2d(\tilde x,\gamma(\tilde x))}\le 8Ne^{2D}\coth^2(\epsilon_3)$$
where $D$ is a the diameter bound on $C(M_Y)^{<\delta}$ and $N$ depends on $D$. By Lemma \ref{diam}, the diameter of each component of $C(M_Y)^{<\delta}$ is bounded by a constant only depending on the the volume of $C(M_Y)$. By Lemma \ref{unif_conn}, $C(M_Y)^{< \delta}$ is connected. Therefore $D$, and hence $N$, only depend on the volume of $C(M_Y)$ and the sum only depends on the volume of $C(M_Y)$. The result then follows from Corollary \ref{thin-area} by further choosing $2\delta < \delta_0$ so that $\Re\cL_\alpha(M_Y) <\delta_0$. \eproof

\section{Vector field at infinity}	
\label{toy-model}
We now return to our flow line $Y_t$ of the gradient vector field $V$ on $\Teich(S)$ and recall that by Theorem \ref{flowfacts} we have $Y_t \to \hat Y \in \Teich(S)$ and $\|V(Y_t)\|_2 = \|\phi_{Y_t}\|_2 \to 0$. The surface $\hat Y$ is a noded surface (with possibly empty nodal set). We can assume that $\alpha$ is one of the nodes and let $c_\alpha(t) = \ell_\alpha(Y_t)/\cL_\alpha(M_{Y_t}).$ Then by Theorem \ref{lemma:c'limit} if $t_n\to \infty$ is a subsequence with $c_\alpha(t_n)\to c$ for some $c \in \CC$ we have $c'_\alpha(t) \to v(c)$ where $v$ is the vector field  
$$v(z) = \frac{1}{4}\left( |z|^4-2z\Re(z^2)-z^2+2z\right).$$

\begin{figure}[htbp] 
   \centering
   \includegraphics[width=4in]{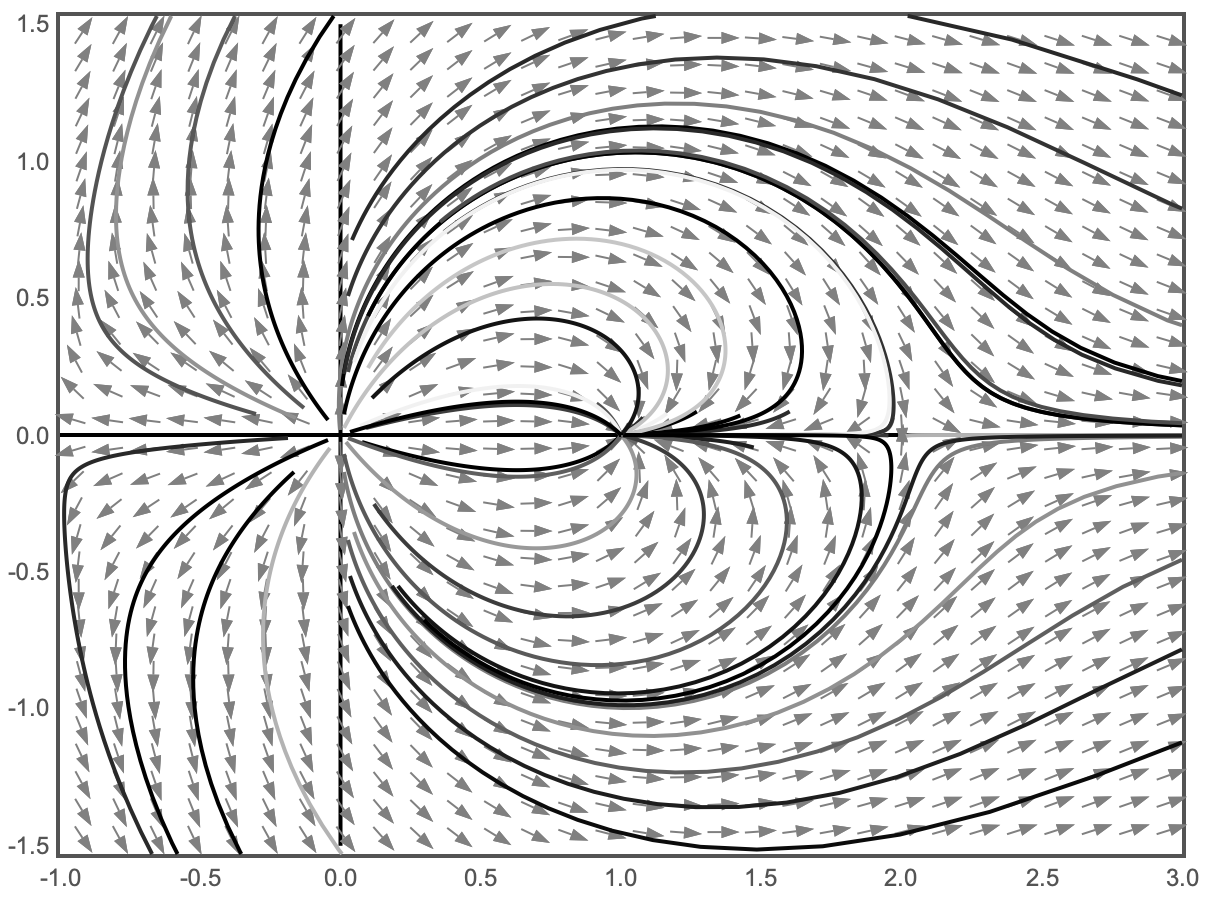} 
   \caption{Vector field $v(z) =\frac{1}{4}\left( |z|^4-2z\Re(z^2)-z^2+2z\right)$  }
   \label{vectorfield}
\end{figure}

\begin{prop}\label{prop:stability}
The critical points of the vector field $v(z)$ are $z=0$ (unstable), $z= 1$ (stable), $z=2$ (saddle), and $z=-1$ (saddle). 
The basin of attraction of $z=1$ is the disk $|z-1| < 1$ and the circle $|z-1| = 1$ consists of two  trajectories  from $0$ to $2$.
\label{lemma:vfield}
\end{prop}

{\bf Proof:} A direct computation shows that $v(-1)=v(0) = v(1) = v(2) = 0$ and that there are no other zeros. The linearization of each zero is non-trivial and gives that each zero of $v$ is as described. Once we show that the circle $|z-1| =1$ consists of two  trajectories  from $0$ to $2$ it will follow that $|z-1|<1$ is the basin of attraction for $z=1$ as there are exactly two trajectories limiting to the unstable fixed point at $z=2$ and therefore all the other trajectories in the disk must limit to $z=1$ (since they can't limit to the unstable fixed point at $z=0$).

If we let $w(z) = (z-1)/|z-1|$ be the radial vector field centered at $z=1$ then the Euclidean inner product of $v$ and $w$ is given by taking the real part of the product of $v$ and the conjugate of $w$ (as functions). Points in the circle $|z-1| = 1$ are of the form $z = 1 + e^{i\theta}$ and on calculates to see that $\langle v(z), w(z) \rangle = 0$ for such points so $v$ is tangent to this circle.

On the other hand if we let $h(z)= 1$ be the constant horizontal vector field then on the circle we have
\begin{equation}\label{horizontal}
\langle v(z), h(z) \rangle = \Re v(1+e^{i\theta}) = \frac34\sin^2\theta(2+\cos\theta).
\end{equation}
This is zero only when $z$ is $0$ or $2$ (and $\theta$ is $\pi$ or $0$) and is positive otherwise. This implies that $z=0$ and $z=2$ are the only critical points on the circle and that the flows lines are from $0$ to $2$ as claimed. \eproof

We will use the next lemma to analyze the limiting behaviour of the path $c_\alpha$.

\begin{lemma}\label{vf_accum}
Let $v$ be a smooth vector field on $\R^n$ and  $\gamma\colon [0,\infty) \to \R^n$ be a smooth path whose image lies in a compact subset. Also assume that for all sequences $t_n$ with $\gamma(t_n)$ converging to some $p\in \R^n$ then $\gamma'(t_n) \to v(p)$.

Let $\mathcal A$ be the accumulation set for $\gamma$. Then $\mathcal A$ is a union of trajectories of the flow of $v$.

If $\cA$ contains distinct points $p$ and $q$ then there exists an $\epsilon>0$ such that if $\delta<\epsilon$ there is a trajectory $\beta \in \cA$ and a $t \in \R$ such that $|\beta(t) - p| = \delta$ and $\langle w, \beta'(t) \rangle \ge 0$ and similarly a trajectory with $\langle w, \beta'(t) \rangle \le 0$. In particular if $p \in \cA$ is an attracting or repelling fixed point of $v$ or $p$ is an isolated point of $\cA$ then $\cA = \{p\}$ and $\gamma(t) \to p$.
\end{lemma}

{\bf Proof:} Since the image of $\gamma$ has compact support we can assume that $v$ has compact support and $|v|$ is bounded. We first show that $|\gamma'(t) - v(\gamma(t))| \to 0$ as $t\to \infty$. If not we can find a sequence $t_n \to \infty$ such that $\gamma(t_n) \to p$ in $M$ but $|\gamma'(t_n) - v(\gamma(t_n))|$ converges to some non-zero $s \in \R$. But, by assumption $\gamma'(t_n) \to v(p)$ which implies $|\gamma'(t_n) - v(\gamma(t_n))| \to 0$, a contradiction.

Therefore $\gamma'(t)$ is bounded. Now let $\gamma_\tau(t) = \gamma(t+\tau)$ and for $p\in \mathcal A$ choose a sequence $\tau_i\to\infty$ with $\gamma_{\tau_i}(0) \to p$. As the $\gamma'_\tau$ are also bounded, by Arzela-Ascoli  (after possibly passing to a subsequence), we have a locally uniform limit $\gamma_i(t) = \gamma_{\tau_i} \to \beta$ and $\gamma'_{\tau_i}(t)$ limits locally uniformly to $v(\beta(t))$. Note that $\beta$ is defined on all of $\R$. Therefore we have
\begin{eqnarray*}
\beta(0) + \int_0^t v(\beta(s))ds & = & \gamma_i(t) - \left(\gamma_i(0) + \int_0^t \gamma'_i(s) ds\right) + \beta(0) + \int_0^t v(\beta(s))ds \\
& =& \lim_{i\to \infty}\left(\gamma_i(t) +\left(\beta(0) -\gamma_i(0)\right) + \int_0^t \left(v(\beta(s)) - \gamma'_i(s)\right) ds\right) \\
& = & \beta(t)
\end{eqnarray*}
which implies that $\beta$ is a trajectory of $v$. Note that the image of $\beta$ will be contained in $\mathcal A$ so we have shown that $\mathcal A$ is a union of trajectories.

For the final statement let $\epsilon = |p-q|/3$. Then for any $\delta>0$ the path $\gamma$ must leave and enter any ball of radius $\delta<\epsilon$ centered at $p$ infinitely often. Therefore we can find $t_n\to \infty$ such that $|\gamma(t_n) - p| = \delta$ and $\langle w, \gamma'(t_n)\rangle \ge 0$. After possibly passing to a subsequence there is a trajectory $\beta \in \cA$ with $\gamma(t_n)\to \beta(t)$ and $\gamma'(t_n) \to \beta'(t)$. Continuity implies that $|\beta(t) - p| = \delta$ and $\langle w, \beta'(t)\rangle \ge 0$ as claimed. The other inequality follows similarly.
\eproof
 
We can now show that $c_\alpha(t)$ converges to one.
 
\begin{prop}
As $t\to \infty$ we have $c_\alpha(t) \to 1$.
\label{lemma:c=1}
\end{prop}
{\bf Proof:}
Let $\cC$ be the accumulation set of $c_\alpha$. As $c_\alpha(t)$ is contained in the disk $|z-1| \le 1$ then we must have all of the trajectories in $\cC$ are contained in this disk. Note that any non-constant trajectory in the disk contains either $0$ or $1$ (and often both)  in its closure so if $\cC$ contains a non-constant trajectory then it contains either $0$ or $1$ and by the second part of Lemma \ref{vf_accum} we then have $\cC = \{0\}$ or $\cC = \{1\}$. In particular, $\cC$ cannot contain any non-constant trajectories and therefore, again by the second part of Lemma \ref{vf_accum}, must be equal to one of the zeros of $v$ in the disk $|z-1| \le 1$. 
We finish the proof by showing that $c_\alpha(t)$ does not converge to $0$ or $2$.

If as $t\to\infty$ we have $c_\alpha(t) \to 0$, 
then  for $f(t) = \Re(\log(c_\alpha(t))$ we have  $f(t) \to -\infty$. But by Theorem \ref{lemma:c'limit} 
$$\lim_{t\rightarrow\infty} \frac{c'_\alpha(t)}{c_\alpha(t)}
  = \frac{1}{2}.$$   
  As $f'(t) = \Re(c_\alpha'(t)/c_\alpha(t))$ it follows that $f(t)$ is increasing for $t$ contradicting thats $\displaystyle\lim_{t\rightarrow\infty}f(t)= -\infty$. Thus $\displaystyle\lim_{t\to\infty} c_\alpha(t) \neq 0$.

If  $c_\alpha(t) \rightarrow 2$ as $t\rightarrow \infty$ then by (1) of Theorem \ref{lemma:c'limit} 
we have 
$$\lim_{n\to\infty} \frac{\ell_\alpha'(t_n)}{\ell_\alpha(t_n)} = \frac32$$
 for all sequences $t_n\to\infty$. On the other hand, as $\ell_\alpha(t)\to 0$ as $t\rightarrow \infty$, we can choose a subsequence such that $\ell_\alpha'(t_n)\leq 0$ and therefore
$$\lim_{n\to\infty} \frac{\ell_\alpha'(t_n)}{\ell_\alpha(t_n)} \le 0.$$
This is a contradiction and therefore $\displaystyle\lim_{t\to\infty}c_\alpha(t) \neq 2$.  It follows that $\displaystyle\lim_{t\to\infty}c_\alpha(t) = 1$. 
\eproof

\section{Main Theorem}
Before we finish the proof of the main theorem, we will need to following bound on the $L^\infty$-norm of a harmonic Beltrami differential in terms of the $L^2$-norm and the derivative of the length of short geodesics. 

\begin{lemma}[{Wolpert, \cite[Lemma 11]{Wolpert:equibounded}}]
There exists an $\epsilon_0$ satisfying the following. Given $\epsilon < \epsilon_0$ there exists $c >0$ such that if  $\tau$ is  the family of geodesics with $\ell_\alpha < \epsilon$ for $\alpha \in \tau$ and $\mu$ is a harmonic Beltrami differential on $Y$ then
$$\|\mu\|_\infty \le c\left(\frac{1}{2}\max_{\alpha\in
\tau} |d(\log\ell_\alpha)(\mu)|+ \|\mu_0\|_2\right)$$
where $\mu_0$ is the component of $\mu$ orthogonal (in the Weil-Petersson inner product) to the span of the of the gradients $\nabla \ell_\alpha$.
\end{lemma}
Wolpert's original bound is in terms of the gradient of the root-length functions $\ell_{\alpha}^{1/2}$. One can translate his statement to the above statement via the chain rule.

The above gives the following immediate corollary;
\begin{corollary}\label{Linfinity_zero}
Let $(N;S)$ be relatively acylindrical and $M_t \in CC(N;S,X)$ be a flowline for $V$, with quadratic differential $\phi_t$. Then 
$$\lim_{t\rightarrow \infty} \|\phi_t\|_\infty = 0$$
\end{corollary}

{\bf Proof:} We have $Y_t \rightarrow Y_\tau$ where $\tau$ is the collection of nodes. Then there is  an $\epsilon_1 >0$ such that if $\ell_\beta(Y_t) < \epsilon_1$ then $\beta\in\tau$.  We apply the above  to $\mu = V_t = -\overline\phi_t/\rho_{Y_t}$ for $\epsilon < \min(\epsilon_0,\epsilon_1)$. As  $\|\mu_0\|_2 \leq \|\mu\|_2$ we have 
$$\|\phi_t\|_\infty \leq c\left(\frac{1}{2}\max_{\alpha\in
\tau} |d(\log\ell_\alpha)(V_t)|+ \|\phi_t\|_2\right).$$
By Lemma \ref{lemma:c=1} $\lim_{t\rightarrow \infty} c_\alpha(t) = c_\alpha = 1$. Thus  by Lemma \ref{lemma:c'limit} 
$$\lim_{t\rightarrow\infty} d(\log\ell_\alpha)(V_t) = \frac{1}{2}(\Re(c_\alpha^2)-1) = 0.$$
Thus as $\displaystyle\lim_{t\rightarrow \infty}\|\phi_t\|_2 = 0$ the result follows.
\eproof

We now have everything in place to prove our main result.
\medskip

\noindent
{\bf Theorem \ref{theorem:main}}
{\em Let $(N;S)$ be relatively acylindrical and $M_t \in CC(N;S,X)$ be a flowline for $V$, then $M_t$ converges to $M_{\rm geod}.$ 
}\newline

{\bf Proof:}
If $Y_t$ is the corresponding flowline of $V$ on $\Teich(S)$ we have that $M_t = M_{Y_t}$. We let $\phi_t$ be the quadratic differential given by the Schwarzian on $Y_t$. Then  by Theorem \ref{flowfacts}  we have that $Y_t \rightarrow \hat Y$ where $\hat Y$ is a possibly noded surface in the Weil-Petersson completion $\overline{\Teich(S)}$. We will  show that 
the set of nodes is empty at the $\hat Y$ is actually contained in $\Teich(S)$.

We let $\sigma= \sigma_{(N;S)}:\Teich(S)\rightarrow \Teich(S)$ be the restriction of the skinning map to $S \subseteq \partial N$. 

By McMullen, the skinning map is contracting in the Teichmuller metric $d_{\Teich}$ with contraction factor $c <1$ depending only on the topology of $(N;S)$ (see \cite[Theorem 6,1, Corollary 6.2]{McMullen:iter}). It follows that there is a unique fixed point of $\sigma$ which we label $Y_{\rm geod}$ since $M_{\rm geod} = M_{Y_{\rm geod}}$. 

For any contraction mapping we can bound the distance from a point to the fixed point in terms of the distance between the point and its first iterate. In particular for any $Y_t$ we have
$$d_{\Teich}(Y_t, Y_{\rm geod}) \le  \frac{d_{\Teich}(Y_t, \sigma(Y_t))}{1-c}.$$
By the Ahlfors-Weill quasiconformal reflection theorem (see \cite[Theorem A]{ahlforsweill}) if $\|\phi_t\|_\infty <1/2$ then 
$$d_{\Teich}(Y_t,\sigma(Y_t)) \leq \frac{1}{2}\log\left(\frac{1+2\|\phi_t\|_\infty} {1-2\|\phi_t\|_\infty}\right).$$
By Corollary \ref{Linfinity_zero} above $\displaystyle\lim_{t\rightarrow \infty}\|\phi_t\| \rightarrow 0$ so we can combine the two inequalities to get $\displaystyle\lim_{t\to\infty}d_{\Teich}(Y_t ,Y_{\rm geod}) =0$. Therefore $M_t \to M_{\rm geod}$ as claimed. \eproof

\bibliography{bib,math}
\bibliographystyle{math}
\end{document}